\documentclass[11pt]{article}

\usepackage[latin1]{inputenc}
\usepackage[T1]{fontenc}
\usepackage{amsmath,amssymb,amsthm}

\DeclareMathOperator{\coker}{coker}
 \DeclareMathOperator{\HH}{H}
 \DeclareMathOperator{\Spec}{Spec}
\DeclareMathOperator{\Hom}{Hom}
\DeclareMathOperator{\Tor}{Tor} \DeclareMathOperator{\Ext}{Ext}
\DeclareMathOperator{\diam}{diam}
 \DeclareMathOperator{\rank}{rank}
\DeclareMathOperator{\Hilb}{Hilb} \DeclareMathOperator{\GradAlg}{GradAlg}

\DeclareMathOperator{\ext}{ext}

\newtheorem{theorem}{Theorem}
\newtheorem{lemma}[theorem]{Lemma}
\newtheorem{corollary}[theorem]{Corollary}
\newtheorem{definition}[theorem]{Definition}
\newtheorem{example}[theorem]{Example}
\newtheorem{proposition}[theorem]{Proposition}

\newtheorem{remark}[theorem]{Remark}

\usepackage{latexsym}            
\usepackage{amssymb}

\newcommand{\pp}{{\mathbb P}}

\newcommand\sF{{\mathcal F}}

\newcommand\sI{{\mathcal I}}

\newcommand\sN{{\mathcal N}}
\newcommand\sO{{\mathcal O}}

\newcommand{\proj}[1]
{ \mathchoice
           { {\mathbb P}^{#1} }
           { {\mathbb P}^{#1} }
           { {\mathbb P}^{#1} }
           { {\mathbb P}^{#1} }
         }

\begin{document}

\title{The Hilbert Scheme of Space Curves of small diameter}

\author{Jan O. Kleppe}


\date{\ }

\maketitle

\vspace*{-0.50in}
\begin{abstract}
\noindent  This paper studies space curves $C$  of degree $d$ and
arithmetic genus $g$, with homogeneous ideal $I$ and Rao module $M =
\HH_{*}^1(\tilde I)$, whose main results deal with curves which satisfy $
{_{0}\!\Ext_R^2}(M ,M ) = 0 $ (e.g. of diameter, $\diam M \leq 2$, which means
that $M$ is non-vanishing in at most two consecutive degrees). For such curves
$C$ we find necessary and sufficient conditions for unobstructedness, and we
compute the dimension of the Hilbert scheme, $\HH(d,g)$, at $(C)$ under the
sufficient conditions. In the diameter one case, the necessary and sufficient
conditions coincide, and the unobstructedness of $C$ turns out to be
equivalent to the vanishing of certain graded Betti numbers of the free graded
minimal resolution of $I$. We give classes of obstructed curves $C$ for which
we partially compute the equations of the singularity of $\HH(d,g)$ at $(C)$.
Moreover by taking suitable deformations we show how to kill certain repeated
direct free factors ("ghost-terms") in the minimal resolution of the ideal of
the general curve. For Buchsbaum curves of diameter at most $2$, we simplify
in this way the minimal resolution further, allowing us to see when a singular
point of $\HH(d,g)$ sits in the intersection of several, or lies in a unique
irreducible component of $\HH(d,g)$. It follows that the graded Betti numbers
mentioned above of a generic curve vanish, and that any irreducible component
of $\HH(d,g)$ is reduced (generically smooth) in the diameter 1 case.

\noindent {\bf AMS Subject Classification.} 14C05, 14H50, 14B10, 14B15, 13D10,
13D02, 13D07, 13C40.


\noindent {\bf Keywords}. Hilbert scheme, space curve, Buchsbaum curve,
unobstructedness, cup-product, graded Betti numbers, ghost term, linkage,
normal module, postulation Hilbert scheme.
 \end{abstract}
 \vspace*{-0.25in}
\thispagestyle{empty}

\section{Introduction and Main Results}

The Hilbert scheme of space curves of degree $d$ and arithmetic genus $g$,
$\HH(d,g)$, has received much attention over the last years after Grothendieck
showed its existence \cite{G}. At so-called special curves it has turned out
that the structure of $\HH(d,g)$ is difficult to describe in detail, and
questions related to irreducibility and number of components, dimension and
smoothness have been hard to solve. For particular classes of space curves,
some results are known. In 1975 Ellingsrud \cite{El} managed to prove that the
open subset of $\HH(d,g)$ of arithmetically Cohen-Macaulay curves (with a
fixed resolution of the sheaf ideal $\sI_C$) is smooth and irreducible, and he
computed the dimension of the corresponding component. A generalization of
this result in the direction of smoothness and dimension was already given in
\cite{K2} (see Theorem~\ref{introth1}$(i)$ below) while the irreducibility was
nicely generalized by Bolondi \cite{B}. Later, Martin-Deschamps and
Perrin gave a stratification $\HH_{\gamma,\rho}$ of $\HH(d,g)$ obtained
by deforming space curves with constant cohomology \cite{MDP1}. Their results
lead immediately to $(iii)$ in the following

%

\begin{theorem} \label{introth1}
  Let $C$ be a curve in $\proj{3}$ of degree $d$ and arithmetic genus $g$, let
  $I = \HH_{*}^0(\sI_C):= \oplus \HH^0(\sI_C(v))$, $M = \HH_{*}^1(\sI_C)$ and
  $E = \HH_{*}^1(\sO_{C})$
  and suppose  
  at least one of the following conditions: \\[-3mm]

  $\ \ \ \ (i) \ \ \ \ \ {_v\!\Hom_R}(I , M) = 0 \ \ \ {\rm for} \ v = 0 \ {\rm
    and} \ v = -4 \ ,$ \\[-4mm]
  
  $\ \ \ \ (ii) \ \ \ \ {_v\!\Hom_R}(M , E) = 0 \ \ {\rm for} \ v = 0 \ {\rm
    and} \ v = -4 \ , \ \ {\rm or}$ \\[-4mm]
  
  $\ \ \ \ (iii) \ \ \ {_0\!\Hom_R}(I , M) = 0 \ , \ \ {_0\!\Hom_R}(M , E) = 0
  \ \ {\rm and} \ \ {_{0}\!\Ext_R^2}(M ,M ) = 0 \ .$  \\[2mm]
  Then $\HH(d,g)$ is smooth at $(C)$, i.e. $C$ is unobstructed. Moreover if
  ${_{0}\!\Ext_R^i}(M,M) = 0$ for $ i \geq 2$, then the dimension of the
  Hilbert scheme at $(C)$ is
\begin{equation*}
  \dim_{(C)}\HH(d,g) = 4d + {_0\!\hom_R}(I,E) + {_{-4}\!\hom_R}(I,M) +
  {_{-4}\!\hom_R}(M , E) \ . 
\end{equation*}
\end{theorem}
We may drop the condition ${_{0}\!\Ext_R^i}(M,M) = 0$ for $ i \geq 2$ in
Theorem~\ref{introth1} by slightly changing the dimension formulas
(cf.\,Theorem~\ref{munobstr} and Remark~\ref{remunobstr}). Moreover we remark
that once we have a minimal resolution of $\sI_C$, we can easily compute $
{_0\!\hom_R}(I,E)$ (as equal to $ \delta^2(0)$ in Definition~\ref{delta})
while the other $_0\!\hom$-dimensions are at least easy to find provided $C$ is
Buchsbaum (Remark~\ref{remunobstr}, \eqref{cor61} and \eqref{cor63}). Another
result from Section 2 is that if a sufficiently general curve $C$ of an
irreducible component $V$ of $\HH_{\gamma,\rho}$ satisfies the vanishing of
the two Hom-groups of Theorem~\ref{introth1}$(iii)$, then $V$ (up to possible
closure in $\HH(d,g)$) is an irreducible component of $\HH(d,g)$
(Proposition~\ref{gen}).

A main goal of this paper is to see when the sufficient conditions of
unobstructedness of Theorem~\ref{introth1} are also {\it necessary}
conditions. Note that it has "classically" been quite hard to prove
obstructedness because one essentially had to compute a neighborhood of $(C)$
in $\HH(d,g)$ to conclude (\cite{Mu}, \cite{K2}, \cite{E}, \cite{GP}). Looking
for another approach to prove obstructedness, we consider in Section 3 the
{\it cup product} and its ``images'' in $ {_0\!\Hom_R}(I,E)$,
${_{-4}\!\Hom_R}(I,M)^{\vee}$ and ${_{-4}\!\Hom_R}(M,E)^{\vee}$ via some
natural maps, close to what Walter and Fløystad do in \cite{W1} and \cite{F}
(see also \cite{MDP3}, \cite{N}). These ``images'' correspond to {\it three}
Yoneda pairings, one of which is the natural morphism
\begin{equation} \label{intro1}
  {_0\!\Hom_R}(I,M) \times {_0\!\Hom_R}(M , E) \longrightarrow  {_0\!\Hom_R}(I ,
  E) \ .
\end{equation}
All three pairings are easy to handle because they are given by taking simple
compositions of homomorphisms, cf.\,Proposition~\ref{cupprod1} and
\ref{cupprod2}. If ${_{0}\!\Ext_R^2}(M,M) = 0$, it turns out that the
non-vanishing of {\it one of} the three pairings is sufficient for
obstructedness. In particular, for a Buchsbaum curve of diameter at most 2, we
can, by using a natural decomposition of $M$, get the non-vanishing of
\eqref{intro1} from the non-vanishing of some of the $\Hom$-groups involved.
More precisely we have (cf. Theorem~\ref{obstr} for a generalization to e.g.
curves with ${_{0}\!\Ext_R^2}(M,M) = 0$ obtained by Liaison Addition)

\begin {theorem} \label{introth2} 
  Let $C$ be a Buchsbaum curve in $\proj{3}$ of diameter at most 2 and let $M
  \cong M_{[c-1]} \oplus M_{[c]}$ be an $R$-module isomorphism where
  $M_{[t]}$, for $t=c-1$ and $c$, is the part of $M = \HH_{*}^1(\sI_C)$
  supported in degree $t$.  Then $C$ is obstructed if one of the following
  conditions hold \\[-3.5mm]

  $ \ \ \ (a) \ \ \ \ \ \ {_0\!\Hom_R}(I ,M_{[t]} ) \neq 0 \ \ \ {\rm and} \ 
  \ \ {_0\!\Hom_R}( M_{[t]} , E) \neq 0 \ , \ \ \ {\rm for} \ t=c \ {\rm or} \ 
  t=c-1 \ , $  \\[-4mm]

$ \ \ \ (b) \ \ \ \ \ {_{-4}\!\Hom_R}(I , M_{[t]}) \neq 0 \ \ \ {\rm and} \ 
\ \ {_0\!\Hom_R}(M_{[t]},E) \neq 0 \ , \ \ \ {\rm for} \ t=c \ {\rm or} \ 
  t=c-1 \ , $  \\[-4mm]

$\ \ \ (c)  \ \ \ \ \ \ \ {_0\!\Hom_R}(I ,  M_{[t]}) \neq 0 \ \ \ {\rm and} 
\ {_{-4}\!\Hom_R}( M_{[t]} , E) \neq 0 \  , \ \ \ {\rm for} \ t=c \ {\rm or}
\  t=c-1 \ . $ \\[-4mm]
\end{theorem}

Buchsbaum curves in $\proj{3}$ are rather well understood by studies of
Migliore and others (cf. \cite{MIG} for a survey of important results as well
as for an introduction to Liaison Addition), and Theorem~\ref{introth2} takes
some care of its obstructedness properties. Note also that since the main
assumption ${_{0}\!\Ext_R^2}(M,M) = 0$ of Section 3 is liaison-invariant,
there may be many more applications of Proposition~\ref{cupprod1} and
\ref{cupprod2}.

Our results in Section 3 also allow an effective calculation of (at least the
degree $2$ terms of) the equations of the singularities of $\HH(d,g)$ at some
curves whose diameter is $2$ or less (as illustrated in Example~\ref{ex1}). To
get equivalent conditions of unobstructedness and a complete picture of the
equations of the singularities of $\HH(d,g)$ more generally, we need a more
general version of the cup product and we certainly need to include their
higher Massey products (Laudal, \cite{L1} and \cite{L2}).

If we reformulate Theorem~\ref{introth1} by logical negation to necessary
conditions of obstructedness (cf. Proposition~\ref{propobstr}) we get
necessary conditions which are quite close (resp. equivalent) to the
sufficient conditions of Theorem~\ref{introth1} in the diameter $2$ case
(resp. in the diameter $1$ case). It is easy to substitute the non-vanishing
of the $\Hom$-groups of Theorem~\ref{introth2} by the non-triviality of
certain graded Betti numbers in the minimal resolution, $$0 \rightarrow
\bigoplus_i R(-i)^{\beta_{3,i}} \rightarrow \bigoplus_i R(-i)^{\beta_{2,i}}
\rightarrow \bigoplus_i R(-i)^{\beta_{1,i}} \rightarrow I \rightarrow 0 \ , $$
of $I$ (cf. Corollary~\ref{mcor}). In the diameter one case, we get the
following main result (cf.\! \cite{MDP1}, pp.\! 185-193 for the case $M \cong
k$).

\begin{theorem} \label{introth3}
  Let $C$ be a curve in $\proj{3}$ whose Hartshorne-Rao module $M \neq 0$ is
  of diameter $1$. 
  Then $C$ is obstructed if and only if
  $$
  \beta_{1,c} \cdot \beta_{2,c+4} \neq 0 \ \ \ {\rm or} \ \ \ \beta_{1,c+4}
  \cdot \beta_{2,c+4} \neq 0 \ \ \ {\rm or} \ \ \ \beta_{1,c} \cdot
  \beta_{2,c}\neq 0 \ . $$
  Moreover if $C$ is {\bf unobstructed} and $M$ is
  $r$-dimensional, then the dimension of the Hilbert scheme $\HH(d,g)$ at
  $(C)$ is
\begin{equation*}
 \dim_{(C)}\HH(d,g) = 4d + {_0\!\hom_R}(I,E) + r(\beta_{1,c+4} + \beta_{2,c}).
\end{equation*}
\end{theorem}

The Hilbert scheme of constant postulation (or the postulation Hilbert
scheme), for which there are various notations, $\GradAlg(H)$,
$\Hilb^H(\proj{3})$ or just $\HH_{\gamma}$ in the literature, has received
much attention recently. We prove

\begin{proposition} \label{introprop}
  In addition to the general assumptions of Theorem~\ref{introth3}, let
  $M_{-4} = 0$. Then
  $$\HH_{\gamma}\ is \ singular \ at \ (C) \ if \ and \ only \ if \ 
  \beta_{1,c+4} \cdot \beta_{2,c+4}\neq 0 \ . $$
  Moreover if $\ \HH_{\gamma}$
  is smooth at $(C)$, then $ \dim_{(C)}\HH_{\gamma} = 4d +
  {_0\!\hom_R}(I,E) + r(\beta_{1,c+4} + \beta_{2,c}- \beta_{1,c}).$
\end {proposition}

In Section 4 we are concerned with curves which admit a generization (i.e. a
deformation to a ``more general curve'') or are
generic in $\HH_{\gamma,\rho}$, $\HH_{\gamma}$ or $\HH(d,g)$. Inspired by
ideas of Martin-Deschamps and Perrin in \cite{MDP1} we prove some results,
telling that we can kill certain repetitions in a minimal resolution
("ghost-terms") of the ideal $I(C)$, under deformation. Hence curves with such
simplified resolutions exist. One result of particular interest is
Theorem~\ref{mainres} which considers the form of a minimal resolution of
$I(C)$ given by a Theorem of Rao, cf.  \eqref{resoluM} and  \eqref{resoluMI}.
We prove  

\begin{theorem} \label{introth5}
 If $C$ is a generic curve of $\ \HH_{\gamma,\rho}$ (or of  $\
  \HH_{\gamma}$ or $\HH(d,g)$), then $C$ admits a minimal free resolution
  of the form 
    $$ 0 \rightarrow L_{4} \stackrel{\sigma \oplus 0}{\longrightarrow} L_3
        \oplus F_1 \rightarrow F_0 \rightarrow I(C) \rightarrow 0 ,
        $$
        where $\sigma: L_4 \rightarrow L_3$ is given by the leftmost map in
        the minimal resolution of the Rao module $M$, cf.  \eqref{resoluM},
        and where $F_1$ and $F_0$ are without repetitions (i.e. without common
        direct free factors).
\end{theorem}

Restricting to general Buchsbaum curves, we prove, under some conditions, 
that $L_4$ and $F_1$, and $L_4$ and $F_0(-4)$, have no common direct free
factor (Proposition~\ref{mainres2}). We get

\begin{corollary} \label{introcor} 
Let $C$ be a curve in $ \proj{3}$ whose Rao module $M \neq
  0$ is of diameter $1$ and concentrated in degree $c$. 
  
  (a) If $C$ is generic in $\HH_{\gamma, \rho}$, then $\HH_{\gamma}$ is 
  smooth at $(C)$. Moreover $C$ is obstructed if and
  only if $\beta_{1,c} \cdot \beta_{2,c+4} \neq 0$. Furthermore if
  $\beta_{1,c} =0$ and $ \beta_{2,c+4}= 0$, then  $C$ is generic in
  $\HH(d,g)$. 
  
  (b) If $C$ is generic in $\HH_{\gamma}$, then $C$ is unobstructed. Indeed
  both
  $\HH(d,g)$ and $\HH_{\gamma}$ are smooth at $(C)$.  In particular every
  irreducible component of $\HH(d,g)$ whose generic curve $C$ satisfies $\diam
  M = 1$ is {\it reduced} (i.e. generically smooth).
\end{corollary}

Moreover we are able to make explicit various generizations of
Buchsbaum curves of diameter at most two, allowing us in many cases to decide
whenever an obstructed curve is contained in a {\it unique\/} component of
$\HH(d,g)$ or not (Proposition~\ref{propfi}). Finally we show that any
Buchsbaum curve whose Hartshorne-Rao module has diameter $2$ or less, admits a
generization in $\HH(d,g)$ to an unobstructed curve, hence belongs to a {\it
  reduced} irreducible component of $\HH(d,g)$. We believe that every
  irreducible component of $\HH(d,g)$ whose generic curve $C$ satisfies $\diam
  M  \leq 2$ is {\it reduced}. 
 
  A first version of this paper (containing Theorem~\ref{munobstr},
  Theorem~\ref{introth3}, Theorem~\ref{introth5}, Corollary~\ref{introcor},
  Proposition~\ref{propfi} and the ``cup product part'' of
  Proposition~\ref{cupprod1} and \ref{cupprod2}, see \cite{K96}, available
  from my home-page) was written in the context of the group "Space Curves" of
  Europroj, and some main results were lectured at its workshop in May 1995,
  at the Emile Borel Center, Paris. Later we have been able to generalize
  several results (e.g. Theorem~\ref{introth2}). The author thanks
  prof. O. A. Laudal at Oslo and prof. G. Bolondi at Bologna for interesting
  discussions on the subject.

\subsection {Notations and terminology} A curve $C$ in $\proj{3}$ is an {\it
  equidimensional, locally Cohen-Macaulay} subscheme of $\proj{}:=\proj{3}$ of
dimension one with sheaf ideal $\sI_C$ and normal sheaf $\sN_C =
\Hom_{\sO_{\proj{}}}(\sI_C,\sO_C)$. If $\sF$ is a coherent
$\sO_{\proj{}}$-Module, we let $\HH^i(\sF) = \HH^i(\proj{},\sF)$,
$\HH_{*}^i(\sF) = \oplus_v \HH^i(\sF(v))$ and $h^i(\sF) = \dim \HH^i(\sF)$,
and we denote by $\chi(\sF) = \Sigma (-1)^i h^i(\sF)$ the Euler-Poincar\'{e}
characteristic. Moreover $M = M(C)$ is the Hartshorne-Rao module
$\HH_{*}^1(\sI_{C})$ or just the Rao module, $E = E(C$) is the module
$\HH_{*}^1(\sO_C)$ and $I=I(C)$ is the homogeneous ideal $\HH_{*}^0(\sI_{C})$ of
$C$. They are graded modules over the polynomial ring $R =
k[X_0,X_1,X_2,X_3]$, where $k$ is supposed to be algebraically closed of
characteristic zero. The postulation $\gamma$ (resp. deficiency $\rho$ and
specialization $\sigma$) of $C$ is the function defined over the integers by
$\gamma(v) = \gamma_C(v) = h^0(\sI_C(v))$ (resp. $\rho(v) = \rho_C(v) =
h^1(\sI_C(v))$ and $\sigma(v) =
\sigma_C(v) = h^{1}(\sO_{C}(v))$). Let  \\[-3mm]

$ \hspace {4cm}  s(C) = \min \{n \arrowvert h^0(\sI_C(n)) \neq 0 \} \ ,  $
\\[-3mm] 

$  \hspace {4cm}  c(C) = \max \{n \arrowvert h^1(\sI_C(n)) \neq 0 \} \ , $
\\[-3mm] 

$ \hspace {4cm} e(C) = \max \{ n \arrowvert h^1(\sO_C(n)) \neq 0 \} \ . $
\\[2mm]
Let $b(C)=\min \{n \arrowvert h^1(\sI_C(n)) \neq 0 \}$ and let $\diam M(C) =
c(C)- b(C)+1 $ be the diameter of $M(C)$ (or of $C$). If $c(C) <
s(C)$ (resp.  $e(C) < b(C)$), we say $C$ has maximal rank (resp.
maximal corank). 
A curve $C$ such that $\mathfrak m \cdot M(C) = 0$, $ \mathfrak m =
(X_0,..,X_3)$, is a Buchsbaum curve. $C$ is {\it unobstructed\/} if the
Hilbert scheme of space curves of degree $d$ and arithmetic genus $g$,
$\HH(d,g)$, is smooth at the corresponding point $(C)= (C \subseteq
\proj{3})$, otherwise $C$ is {\it obstructed}. The open part of $\HH(d,g)$ of
{\it smooth connected} space curves is denoted by $\HH(d,g)_S$, while
$\HH_{\gamma,\rho} = \HH(d,g)_{\gamma,\rho}$ (resp. $\HH_{\gamma}$, resp.
$\HH_{\gamma,M})$ denotes the subscheme of $\HH(d,g)$ of curves with constant
cohomology, i.e. $\gamma_C$ and $\rho_C$ do not vary with $C$, (resp. constant
postulation $\gamma$, resp. constant postulation $\gamma$ and constant Rao
module $M$), cf. \cite{MDP1} for an introduction. The curve in a sufficiently
small open irreducible subset of $\HH(d,g)$ (small enough to satisfy all the
openness properties which we want to pose) is called a {\it generic} curve of
$\HH(d,g)$, and accordingly, if we state that a generic curve has a certain
property, then there is an non-empty open irreducible 
subset of $\HH(d,g)$ of curves having this property.  A {\it generization}
$C' \subseteq \proj{3}$ of $C \subseteq \proj{3}$ in 
$\HH(d,g)$ is a generic curve of some irreducible subset of $\HH(d,g)$
containing $(C)$.

For any graded R-module $N$, we have the right derived functors $\HH_{\mathfrak
  m}^i(N)$ and ${_v\!\Ext_{\mathfrak m}^i}(N , -)$ of $\Gamma_{\mathfrak m}(N)
= \oplus_v \ker(N_v \rightarrow \Gamma( \proj{}, \tilde N (v)))$ and
$\Gamma_{\mathfrak m}(\Hom_R(N , -))_v$ respectively (cf. \cite{SGA2}, exp.
  VI). 
  We use small letters for the $k$-dimension and subscript $v$ for the
  homogeneous part of degree $v$, e.g.  ${_v\!\ext_{\mathfrak m}^i}(N_1,N_2) =
  \dim {_v\!\Ext_{\mathfrak m}^i}(N_1,N_2)$.

\section {Preliminaries. Sufficient conditions for unobstructedness.}

In this section we recall the main Theorem on unobstructedness of space curves
of this paper (Theorem~\ref{introth1} or Theorem~\ref{munobstr}).
Theorem~\ref{munobstr} is not entirely new. Indeed $(i)$ and $(i')$ were
proved in \cite{K2} under the assumption ``$C$ generically a complete
intersection'' (combining \cite{K79}, Rem. 3.7 and \cite{K3}, (4.10.1) will
lead to a proof), while the $(iii)$ and $(iii')$ part is a rather
straightforward consequence of a theorem of Martin-Deschamps and Perrin which
appeared in \cite{MDP1}. However, $(ii)$ and $(ii')$ seem new, even though at
least $(ii)$ is easily deduced from $(i)$ by linkage. Indeed linkage preserves
unobstructedness also in the non arithmetically Cohen-Macaulay (ACM) case
provided we link carefully (Proposition~\ref{mlink}). We will include a proof
of Theorem~\ref{munobstr}, also because we need the arguments (e.g. the
technical tools and the exact sequences which appear) later.

Let $N$, $N_1$ and $N_2$ be graded $R$-modules of finite type. Then recall
that the right derived functors $_v\!\Ext_{\mathfrak m}^i (N,-)$ of
$_v\!\HH_{\mathfrak m}^0 (\Hom_R (N,-))$ are equipped with a spectral sequence
(\cite{SGA2}, exp. VI)
\begin{equation}\label{mspect}
E_2^{p,q} = {_v\!\Ext_R^p }(N_1 , \HH_{\mathfrak m}^q (N_2 )) \Rightarrow
{_v\!\Ext_{\mathfrak m}^{p+q} }(N_1, N_2 )
\end{equation}
($\Rightarrow$ means ``converging to'') and a duality isomorphism (\cite{K5},
Thm. 1.1);
\begin{equation}\label{mduality}
{_v\!\Ext^i_{\mathfrak m}} (N_2 ,N_1 ) \cong \hbox{$_{-v-4}
\Ext_R^{4-i} (N_1 ,N_2 )^\vee$}
\end{equation}
where $(-)^\vee = \Hom_k( - , k)$, which generalizes the Gorenstein duality
$_{v}\!\HH_{\mathfrak{m}}^{i}(M) \simeq \ 
_{-v}\!\Ext_R^{4-i}(M,R(-4))^{\vee}.$ These groups fit into a long exact
sequence (\cite{SGA2}, exp.\ VI)
\begin{equation}\label{seq}
\rightarrow {_v\!\Ext_{\mathfrak m}^i }(N_1 ,N_2 ) \rightarrow {_v\!\Ext_R^i }
(N_1 ,N_2 ) \rightarrow \Ext_{{\mathcal O}_{\proj{}}}^i (\tilde N_1
,\tilde N_2 (v)) \rightarrow  {_v\!\Ext_{\mathfrak m}^{i+1} }(N_1 ,N_2 )
\rightarrow
\end{equation}
which in particular relates the deformation theory of $(C \subseteq
\proj{3})$, described by $\HH^{i-1}({\sN_C}) \cong \Ext_{{\mathcal
    O}_{\proj{}}}^i (\tilde I ,\tilde I (v))$ for $i = 1,2$ (cf. \cite{K2},
Rem. 2.2.6 for a proof of this isomorphism), to the deformation theory of the
homogeneous ideal $I = I(C)$, described by ${_0\!\Ext_R^i } (I ,I )$, in the
following exact sequence
\begin{equation}\label{mseq}
{_v\!\Ext_R^1 }(I,I ) \hookrightarrow  \HH^0({\sN_C}(v)) \rightarrow
{_v\!\Ext_{\mathfrak m}^2}(I,I ) \stackrel{\alpha}{\longrightarrow}
{_v\!\Ext_R^2 }(I,I ) \rightarrow 
\HH^1({\sN_C}(v)) \rightarrow {_v\!\Ext_{\mathfrak m}^3 }(I ,I ) \rightarrow 0
\ .
\end{equation}
Let $M= \HH_{\mathfrak m}^2(I)$. In this situation C. Walter proved that the
map $\alpha: {_v\!\Ext_{\mathfrak m}^2 }(I ,I ) \cong {_v\!\Hom_R}(I
,\HH_{\mathfrak m}^2(I) ) \rightarrow {_v\!\Ext_R^2 }(I,I )$ of \eqref{mseq}
factorizes via $ {_v\!\Ext_R^2 }(M , M)$ in a natural way (\cite{W2}, Thm.
2.3), the factorization is in fact given by a certain edge homomorphism of the
spectral sequence \eqref{mspect} with $N_1 = M$, $N_2=I$ and $p+q = 4$, cf.
\eqref{thm63}, \eqref{thm64} and \eqref{thm65} where the factorization of this
map occurs. Fløystad furthered the study of $\alpha$ in \cite{F}. Also in
\cite{MDP3}, (see \cite{MDP3} Sect.\! 0.e and Sect.\! 3), they need to
understand $\alpha$ properly to make their calculations.

To compute the dimension of the components of $\HH(d,g)$, we have found it
convenient to introduce the following invariant, defined in terms of the
graded Betti numbers of a minimal resolution of the homogeneous ideal $I$ of
$C$:
\begin{equation} \label{resolu}
0 \rightarrow \bigoplus_i R(-i)^{\beta_{3,i}} \rightarrow \bigoplus_i
R(-i)^{\beta_{2,i}} \rightarrow \bigoplus_i R(-i)^{\beta_{1,i}} \rightarrow I
\rightarrow 0 
\end{equation}

\begin{definition} \label{delta} 
  If $C$ is a curve in $\proj{3}$, we let
\begin{equation*}
  \delta^{j}(v) = \sum _{i} \beta_{1,i} \cdot h^j(\sI_C(i+v)) - \sum _{i}
  \beta_{2,i} \cdot h^j(\sI_C(i+v))+  \sum _{i} \beta_{3,i} \cdot
  h^j(\sI_C(i+v))   
\end{equation*}
\end{definition}

\begin{lemma}\label{eulerI}
  Let $C$ be any curve of degree $d$ in $\proj{3}$. Then the following
  expressions are equal
\begin{equation*}
 {_0\!\ext_R^1 }(I,I )- {_0\!\ext_R^2 }(I,I ) = 1 - \delta^0(0) = 4d +
\delta^2(0) - \delta^1(0) = 1 + \delta^2(-4) - \delta^1(-4)
\end{equation*}
\end{lemma}

\begin{remark} \label{remeulerI} 
  Those familiar with results and notations of \cite{MDP1} will recognize $1-
  \delta^{0}(0)$ as $\delta_{\gamma}$ and $\delta^{1}(-4)$ as
  $\epsilon_{\gamma ,\delta}$ in their terminology. By Lemma~\ref{eulerI} it
  follows that the dimension of the Hilbert scheme $\HH_{\gamma,M}$ of
  constant postulation and Rao module, which they show is $\delta_{\gamma} +
  \epsilon_{\gamma,\delta} -{_0\!\hom}(M,M )$ (Thm. 3.8, page 171), is
  also equal to $1+\delta^{2}(-4)- {_0\!\hom}(M,M)$.
\end{remark}

\begin {proof}  To see the equality to the left, we apply
  ${_v\!\Hom_R}(-, I)$ to the resolution \eqref{resolu}. Since $\Hom_R(I , I)
  \cong R $ and since the alternating sum of the dimension of the terms in a
  complex equals the alternating sum of the dimension of its homology groups,
  we get
\begin{equation} \label{lem21}
  {dim R_v} - {_v\!\ext_R}^1(I , I) +
{_v\!\ext_R}^2(I , I) = \delta^0(v) \ \  , \ \ \ \   v \in Z
\end{equation}
If $v = 0$ we get the equality of Lemma~\ref{eulerI} to the left. The
equality in the middle follows from \cite{K2}, Lemma 2.2.11. We will, however,
indicate how we can prove this and the right hand equality by using 
\eqref{mduality} and \eqref{seq}. Indeed by \eqref{mduality},
${_v\!\ext_m^{4-i}}(I , I) = {_{-v-4}\!\ext_R^i}(I , I)$. Hence
\begin{equation}  \label{lem22}
  {_v\!\ext_m^2}(I , I) - {_v\!\ext_ m^3}(I , I) + dim R_{-v-4} =
  \delta^0(-v-4) \ \ ,\ \ \ v \in Z  
\end{equation}
by \eqref{lem21}. Combining \eqref{lem21} and \eqref{lem22} with the exact
sequence  \eqref{mseq}, we get
\begin{equation}  \label{lem23}
{ v+3  \choose 3} ~ - ~  \chi (\sN _{C} (v)) ~ =  ~ \delta ^{0} (v) ~ - ~
\delta ^{0} (-v-4) ~ ~ ~ , ~ ~ ~ v ~  \in ~ Z 
\end{equation}
because $\dim R_v - \dim R_{-v-4} = { v+3 \choose 3}$. Therefore it suffices
to prove
\begin{equation} \label{lem24}
\delta ^{0}(-v-4) =   \delta ^{1}(v) - \delta ^{2}(v) ~ ~ , ~ ~ ~ v
~\ge ~ -4 ~ ~ ~ ~ ~ 
\end{equation}
Indeed using \eqref{lem23} and \eqref{lem24} for $v = 0$ we get the equality
of Lemma~\ref{eulerI} in the middle because $\chi(\sN_ C) = 4d$ holds for {\it
  any\/} curve (cf. Remark~\ref{remeulerI2}) while \eqref{lem24} for $v = -4$
takes care of the last equality appearing in  Lemma~\ref{eulerI}. 

To prove \eqref{lem24} we use the spectral sequence \eqref{mspect} together
with \eqref{lem22}. Letting $M = \HH_{\mathfrak m}^2(I)$ and $E =
\HH_{\mathfrak m}^3(I)$ we get ${_v\!\Ext_{\mathfrak m}^2}(I , I) \cong
{_v\!\Hom_R}(I, M)$ and ${_v\!\Ext_R^2}(I , E) \cong{_v\!\Ext_{\mathfrak
    m}^5}(I , I) = 0$ and an exact sequence
\begin{equation}  \label{lem25} 
 {_v\!\Ext_R^1}(I, M) \hookrightarrow {_v\!\Ext_{\mathfrak
  m}^3}(I,I) \rightarrow  {_v\!\Hom_R}(I , E)  \rightarrow   {_v\!\Ext
  _R^2}(I,M) \rightarrow  {_v\!\Ext_{\mathfrak m}^4}(I,I) \twoheadrightarrow  
  {_v\!\Ext_R^1}(I , E)
\end{equation}
where we have used that $v \geq -4$ implies ${_v\!\Hom}(I, \HH_{\mathfrak
  m}^4(I)) = 0$ since $ \HH_{\mathfrak m}^4(I)= \HH_{\mathfrak m}^4(R)$. As
argued for \eqref{lem21}, applying ${_v\!\Hom}(-, M)$ (resp. ${_v\!\Hom}(-,
E)$) to the resolution \eqref{resolu}, we get
\begin{equation}  \label{lem26} 
\delta ^{1} (v)  =  \sum_{ i = 0 }^{2} (-1) ^{i} ~ {_v\!\ext ^i} (I
 , M) ~ ~ ~ , ~ ~ ~( ~ {\rm resp.} ~ ~ \delta ^{2} (v)  =  \sum_{ i = 0
}^{2} (-1) ^{i} ~ {_v\!\ext ^i} (I  ,  E) ~ ) 
\end{equation}
So $\delta^1(v)-\delta^2(v)$ equals $\sum_{ i = 2 }^{4} (-1)^{i} \,
{_v\!\ext_{\mathfrak m}^i} (I , I)$ by \eqref{lem25}, and since
${_v\!\Ext_{\mathfrak m}^4}(I , I) \cong {_{-v-4}\!\Hom}(I,I)^{\vee} \cong
R_{-v-4}^{\vee}$ we get \eqref{lem24} from \eqref{lem22}, and the proof of
Lemma~\ref{eulerI} is complete.
\end{proof}

\begin{remark} \label{remeulerI2} In \cite{K2}, Lemma 2.2.11 we proved
  $\chi(\sN_C(v))=2dv + 4d$ for {\it any} curve and any integer $v$ by
  computing $\delta^0(v)$ for $v >> 0$. Indeed using the definition of
  $\delta^0(v)$, the sequence $ 0 \rightarrow \sI_C \rightarrow \sO_{\proj{}}
  \rightarrow \sO_C \rightarrow 0 $ and $\sum_j \ (-1)^j \sum_i \  i \cdot
  \beta_{j,i} = 0$, we get by applying Riemann-Roch to $\chi(\sO_C(i+v))$,

\begin{equation}  \label{lem27} 
 \delta^0(v) = \sum_j \sum_i (-1)^j \beta_{j,i} \cdot \chi(\sO_{\proj{}}(i+
    v)) - (dv + 1 -  g)  ~ ~ , ~ ~ ~  v >> 0
\end{equation}
while duality on $\proj{}$ and \eqref{resolu} show that the double sum of
\eqref{lem27} equals $-\chi(\sI_C(-v-4))= { v+3 \choose 3}+ \chi
(\sO_{C}(-v-4))$. We get
$\chi(\sN_C(v)) = 2dv + 4d$ by combining with \eqref{lem23}.
\end{remark} 

\begin{proposition}\label{mlink}
  Let $C$ and $C'$ be curves in $\proj{3}$ which are linked (algebraically)
by a complete intersection of two surfaces of degrees $f$ and $g$. If
\begin{equation*}
\HH^1(\sI_C(v)) = 0 \ \ \ {\rm for} \ v= f , g , f - 4 \ {\rm and} \  g - 4,
\end{equation*}
then $C$ is unobstructed if and only if $C'$ is unobstructed.
\end{proposition}

One may find a proof in \cite{K3}, Prop. 3.2. Proposition~\ref{mlink} allows
us to complete the proof of the following main result on unobstructedness. It
applies mostly to curves of small diameter, see also Mir\'o-Roig's criterion
for unobstructedness of Buchsbaum curves of diameter at most 2 (\cite{Mir}).

\begin{theorem} \label{munobstr}
  If $C$ is any curve in $\proj{3}$ of degree $d$ and
  arithmetic genus $g$, satisfying (at least) one of the following conditions:
  \\[-3.5mm] 

  $\ \ \ \ (i) \ \ \ \ \ {_v\!\Hom_R}(I , M) = 0 \ \ \ {\rm for} \ v = 0 \ {\rm
    and} \ v = -4$ \\[-4mm]
  
  $\ \ \ \ (ii) \ \ \ \ {_v\!\Hom_R}(M , E) = 0 \ \ {\rm for} \ v = 0 \ {\rm
    and} \ v = -4$  \\[-4mm]
  
  $\ \ \ \ (iii) \ \ \ {_0\!\Hom_R}(I , M) = 0 \ ,\ {_0\!\Hom_R}(M , E) = 0 \
  \  {\rm and} \ \ {_0\!\Ext_R^2}(M , M) = 0 \ ,$ \\[2mm]
  then $C$ is unobstructed. Moreover, in each case, the dimension of the
  Hilbert scheme $H(d,g)$ at $(C \subseteq \proj{3})$ is given by  \\[-4mm]
  
  $\ \ \ (i') \ \ \ \  \dim_{(C)} H(d,g) = 4d + \delta^2(0) - \delta^1(0) \ 
  ,\ \ {\rm provided} \ (i) \ {\rm holds} \ ,$  \\[-4mm]

$ \ \ \ (ii')  \ \ \  \dim_{(C)} H(d,g) = 4d + \delta^2(0) - \delta^1(0) +
{_{-4}\!\hom_R}(I,M) + {_0\!\hom_R}(I,M) - {_0\!\ext_R^2}(M,M), \\
{\hspace{2.2 cm}} {\rm provided} \ \ (ii) \ \ {\rm holds},$  \\[-4mm]

$ \ \ \ (iii')   \ \  \dim_{(C)} H(d,g) = 4d + \delta^2(0) - \delta^1(0)+
   {_{-4}\!\hom_R}(I,M) , \ \ {\rm provided} \  (iii) \ {\rm holds}. $
\end{theorem}

\begin{proof} $(i)$ Let $A=R/I$ and let $Def_{I}$ (resp. $Def_{A}$) be the
  deformation functor of deforming the homogeneous ideal $I$ as a graded
  $R$-module (resp.\! $A$ as a graded quotient of $R$), defined on the
  category of local Artin $k$-algebras with residue field $k$. Let $\Hilb_{C}$
  be the corresponding deformation functor of $C \subseteq \proj{3}$ (i.e the
  local Hilbert functor at $C$) defined on the same category. To see that $C$
  is unobstructed we just need, thanks to the duality \eqref{mduality}, to
  interpret the exact sequence \eqref{mseq} in terms of deformation theory.
  Recalling that $\sO_{C,x}$, $x\in C$ is unobstructed since $\sI_{C,x}$ has
  projective dimension one (cf. \cite{El}), we get that $\HH^1({\sN_C})$
  contains all obstructions of deforming $C \subseteq \proj{3}$. By
  \eqref{mspect} and \eqref{mduality};
\begin{equation} \label{thm61}
  {_0\!\Ext_{\mathfrak m}^2}(I , I) \cong {_0\!\Hom}(I , M) \ ,\ \
    {\rm  and}\ \  {_0\!\Ext_R^2}(I, I) \cong{_{-4}\!\Ext_{\mathfrak
    m}^2}(I,I)^{\vee}  \cong {_{-4}\!\Hom}(I, M)^{\vee} \ . 
\end{equation} 
Using the vanishing of the first group of \eqref{thm61}, we get $Def_{I} \cong
\Hilb_{C} $ since \eqref{mseq} shows that their tangent spaces are isomorphic
and since we have an injection of their obstruction spaces (similar to the
proof of $ Def_{A} \cong \Hilb_{C}$ in \cite{K79}, Rem.\! 3.7, where the
{\it former} functor must be isomorphic to the local Hilbert functor of
constant postulation of $C$ because it deforms the graded quotient $A$ flatly,
i.e. has constant Hilbert function), cf. \cite{K2}, Thm.\! 2.2.1 and
\cite{W1}, Thm.\! 2.3 where Walter manages to get rid of the ``generically
complete intersection'' assumption of \cite{K2}, § 2.2 by the argument in the
line before \eqref{thm61} (see also \cite{F}, Prop.\! 3.13 or \cite{MDP1},
VIII, for their tangent spaces). Now $Def_{I}$ is smooth because
${_0\!\Ext_R^2}(I,I)$ vanishes by \eqref{thm61}.
This proves $(i)$, and then $(i')$ follows at once from Lemma~\ref{eulerI}.

$(iii)$ One may deduce the unobstructedness of $C$ from results in \cite{MDP1}
by combining Thm.\! 1.5, page 135 with their tangent space descriptions, pp.\!
155-166. However, since we need the basic exact sequences below later (for
which we have no complete reference), we give a new proof which also leads to
another result (Proposition~\ref{gen}$(b)$). Indeed for any curve we claim
there is an exact sequence:
\begin{equation} \label{thm62}
0 \rightarrow  T _{ \gamma  , \rho}  \rightarrow  {_0\!\Ext_R^1} (I , I)
\stackrel{\beta}\longrightarrow  {_0\!\Hom_R} (M , E) \rightarrow
{_0\!\Ext_R^2} (M ,M) \rightarrow   {_0\!\Ext_R^2}(I, I) \rightarrow 
\end{equation}
where $T_{\gamma , \rho}$ is the tangent space of the Hilbert scheme of
constant cohomology $\HH_{\gamma ,\rho}$ at $(C)$. To prove it we use the
spectral sequence \eqref{mspect} and the duality \eqref{mduality} {\it
  twice\/} (Walter's idea mainly, to see the factorization of $\alpha$ via
${_0\!\Ext_R^2} (M ,M)$ in \eqref{mseq}), to get an isomorphism, resp. a
surjection
\begin{equation}  \label{thm63}
  {_0\!\Ext_R^2}(I,I) \cong  {_{-4}\!\Ext_{\mathfrak m}^2}(I, I)^{\vee} \cong 
      {_{-4}\!\Hom}(I, M)^{\vee}  \cong {_0\!\Ext_{\mathfrak m}^4}(M, I) 
\end{equation}
\begin{equation}  \label{thm64}
  \beta_1: {_0\!\Ext_R^1}(I,I) \cong {_{-4}\!\Ext_{\mathfrak m}^3}(I,I)^{\vee}
    \twoheadrightarrow {_{-4}\!\Ext_R^1}(I,M)^{\vee} \cong {_0\!\Ext_{\mathfrak
    m}^3}(M , I)
\end{equation}
Now replacing $I$ by $M$ as the first variable in \eqref{lem25} or using
\eqref{mspect} directly, we get
\begin{equation}  \label{thm65}
  0  \rightarrow   {_0\!\Ext_R^1}(M,M) \rightarrow   {_0\!\Ext_{\mathfrak
    m}^3}(M,I) \stackrel{\beta_2}\longrightarrow  {_0\!\Hom}(M,E) \rightarrow
    {_0\!\Ext_R^2}(M,M) \rightarrow  {_0\!\Ext_{\mathfrak m}^4}(M,I)
    \rightarrow 
\end{equation}
which combined with \eqref{thm63} and \eqref{thm64} yield \eqref{thm62}
because the composition $\beta$ of $\beta_1$ (arising from duality used twice)
and $\beta_2$ must be the natural one, i.e. the one which sends an extension
of ${_0\!\Ext^1}(I,I)$ (i.e. a short exact sequence) onto the corresponding
connecting homomorphism $M = \HH_{\mathfrak m}^2(I) \rightarrow E =
\HH_{\mathfrak m}^3(I)$. And we get the claim by \cite{MDP1}, Prop.\! 2.1, page
157, which implies  $\ker \beta = T_{\gamma, \rho}$.

To see that $C$ is unobstructed, we get by \eqref{thm62} and the vanishing of
${_0\!\Hom_R}(M , E)$ an isomorphism between the local Hilbert functor of
constant cohomology at $C$ and $Def_{I}$. The latter functor $Def_{I}$ is
isomorphic to $\Hilb_{(C)}$ because ${_0\!\Hom_R}(I , M) = 0$ (cf. the proof
of $(i)$), while the former functor is smooth because ${_0\!\Ext^2}(M , M)$
contains in a natural way the obstructions of deforming a curve in
$\HH_{\gamma,\rho}$ (cf. \cite{MDP1}, Thm. 1.5, page 135). This leads easily to
the conclusion of $(iii)$. Moreover note that we now get $(iii')$ from
Lemma~\ref{eulerI} because $h^0(\sN_{C}) ={_0\!\ext_R^1}(I,I)$ and
${_0\!\ext_R^2}(I, I)= {_{-4}\!\hom_R}(I , M)$.
 
$(ii)$ The unobstructedness of $C$ follows from Proposition~\ref{mlink}. Indeed
if we take a complete intersection $Y \supseteq C$ of two surfaces of degrees 
$f$ and $g$ such that the {\it conditions\/} of Proposition~\ref{mlink} hold
(such $Y$ exists), then the corresponding linked curve $C'$ satisfies
\begin{equation} \label{thm66}
  {_v\!\Hom_R}(I(C') , M(C')) \cong {_v\!\Hom_R}(M(C),E(C)) \ \ {\rm for} \ v
  = 0 \ {\rm and} \ v = -4
\end{equation}
because $M(C')$ (resp. $I(C')/I(Y)$) is the dual of $M(C)(f+g-4)$ (resp.
$E(C)(f+g-4)$) and ${_v\!\Hom_R}(I(Y),M(C'))= 0$ for $v=0,-4$ by assumption.
Hence we conclude by Proposition~\ref{mlink} and Theorem~\ref{munobstr}$(i)$.
It remains to prove the dimension formula in $(ii')$. For this we claim that
the image of the map $\alpha : {_0\!\Ext_{\mathfrak m}^2}(I,I) \cong
{_0\!\Hom_R}(I,M) \rightarrow {_0\!\Ext_R^{2}}(I,I)$ which appears in
\eqref{mseq} for $v = 0$, is isomorphic to ${_0\!\Ext_R^2}(M,M)$. Indeed
$\alpha$ factorizes via ${_0\!\Ext_R^2}(M,M)$ in a natural way, and the
factorization is given by a certain map of \eqref{thm62}.  Now ${_v\!\Hom_R}(M
, E) = 0$ for $v = 0$ and $-4$ implies that the maps ${_v\!\Ext_R^2}(M, M)
\rightarrow {_v\!\Ext_R^2}(I, I)$ of \eqref{thm62} are injective for $v = 0$
and $v = -4$. Dualizing one of them (the map for $v = -4$) we get a surjective
composition;
\begin{equation} \label{thm67}
  {_0\!\Hom_R}(I , M) \cong {_{-4}\!\Ext_R^2}(I,I)^{\vee} \rightarrow
        {_{-4}\!\Ext_R^2}(M , M)^{\vee} \cong {_0\!\Ext_R^2}(M, M)
\end{equation}
 which composed with the other injective map above is precisely
$\alpha$. This proves the claim. Now by \eqref{mseq} and the proven claim;
$$h^0(\sN_C) = {_0\!\ext_R^1}(I , I) + \dim \ker \alpha = {_0\!\ext_R^1}(I ,
I) + {_0\!\hom_R}(I , M) - {_0\!\ext_R^2}(M, M)$$
and we get the dimension formula by Lemma~\ref{eulerI} and we are done.
\end{proof}

\begin{remark}  \label{remunobstr}
  Putting the arguments in the text at \eqref{lem25} and \eqref{lem26}
  together (and use that $ {_0\!\Ext_{\mathfrak m}^4}(I,I)=0$), we get
\begin{equation} \label{delta2}  
\delta^{2}(0) = {_0\!\hom_R}(I , E)\ . 
\end{equation}
Moreover if
\begin{equation} \label{exti}
 {_0\!\Ext_R^i}(M , M) = 0 \ \ {\rm for } \ \  2 \leq i \leq 4 \ .
\end{equation}
then we may put the different expressions of $\dim_{(C)}\HH(d,g)$ of
Theorem~\ref{munobstr} in one common formula;
\begin{equation}
\dim_{(C)}\HH(d,g) = 4d + \delta^{2}(0) + {_{-4}\!\hom_R}(I,M) +
{_{-4}\!\hom_R}(M, E).
\end{equation}
Indeed $ {_0\!\Ext_R^2}(I, M) \cong {_{-4}\!\Ext_{\mathfrak m}^2}(M, I)^{\vee}
\cong {_{-4}\!\Hom}(M , M)^{\vee} \cong {_0\!\Ext_R^4}(M, M) = 0, $ and we
have $${_0\!\Ext_R^1}(I,M) \cong {_{-4}\!\Ext_{\mathfrak m}^3}(M,I)^{\vee}
\cong {_{-4}\!\Hom}(M,E)^{\vee}$$
because ${_{-4}\Ext_R^i}(M,M) \cong
{_0\Ext_R^{4-i}}(M,M)^{\vee} = 0$ for $i = 1,2$. Hence $
{_{0}\!\ext_{\mathfrak m}^3}(I,I) = {_0\!\hom_R}(I ,E) + {_{-4}\!\hom}(M,E)$
by \eqref{lem25}. Using \eqref{thm61} and that $\alpha=0$ in \eqref{mseq} for
$v=0$, we get
\begin{equation}  \label{thm51} 
h^1(\sN_C) =  \delta^2(0) + {_{-4}\!\hom_R}(I,M) + {_{-4}\!\hom}(M,E) \ ,
\end{equation}
and we conclude easily.
\end{remark}

Using \eqref{thm51}, we can generalize the vanishing result of $\HH^1(\sN_C)$
appearing in \cite{K3}, Cor.\! 4.12, to

 \begin{corollary} \label{H1N}
   Let $C$ be any curve in $\proj{3}$, let $\diam M \leq 2$ and suppose $e(C)
   < s(C)$. If $\diam M \neq 0$, suppose also $e(C) \leq c+1-\diam M$ and
   $c(C) \leq s(C)$. Then $$\HH^1(\sN_C)=0 \ .$$
\end{corollary}

\begin{proof}
  Since $e(C) < s(C)$, we get $\delta^{2}(0) = 0$ by the definition of
  $\delta^{2}(0)$. Moreover suppose $C$ is not ACM. Then $c(C) \leq s(C)$ and
  \eqref{resolu} imply ${_0\!\ext_R^1}(I,M) = 0$. Finally, since we have $
  \max \{ i \, \arrowvert \, \beta_{1,i} \neq 0 \} \leq \max \{\,
  c(C)+2,e(C)+3\, \}$ by Castelnuovo-Mumford regularity, we get
  ${_{-4}\!\hom_R}(I,M)=0$ by \eqref{resolu} and we conclude by
  Remark~\ref{remunobstr}.
\end{proof}  

Hence curves of $\diam M \leq 2$ whose minimal resolution \eqref{resolu} is
``close enough'' to being linear satisfy $\HH^1(\sN_C)=0$. Indeed
$\HH^1(\sN_C)=0$ for any curve of diameter one or two (resp. diameter zero)
whose Betti numbers satisfy $\beta_{2,i} = 0$ for $i > \min \{c+5-\diam
M,s+3\}$, $\beta_{1,i} = 0$ for $i < c$ (resp. $\beta_{2,i} = 0$ for $i >
s+3$). Thus Corollary~\ref{H1N} generalizes \cite{MN}, Prop.\! 6.1. 

\begin{remark} \label{remmunobstr}
  \eqref{thm61}, \eqref{thm62}, \eqref{thm63}, \eqref{thm64} and
  \eqref{thm65} are valid for {\it any} curve in $\proj{3}$. Moreover if
  $M_{-4} = 0$, we get ${_0\!\Hom(M,\HH_{\mathfrak m}^4}(I)) \cong
  {_0\!\Ext_{\mathfrak m}^4}(M ,R) = 0$ since
  $\HH_{\mathfrak m}^4(I) \cong \HH_{\mathfrak m}^4(R)$ and one may see
that the spectral sequence which converges to ${_0\!\Ext_{\mathfrak m}^4}(M ,
I)$ (cf.  \eqref{thm63}, \eqref{thm64} and \eqref{thm65}) consists of at most
two non-vanishing terms. Hence we can continue the exact sequences
\eqref{thm65} and \eqref{thm62} to the right with    
    $${_0\!\Ext_{\mathfrak m}^4}(M , I) \cong {_0\!\Ext^2}(I, I) \rightarrow
    {_0\!\Ext_R^1}(M , E) \rightarrow\ {_0\!\Ext_R^3}(M , M) $$
 \end{remark}
 
 The proof of Theorem~\ref{munobstr} implies also the following result (see
 $(i)$, mainly the argument from \cite{K79}, Rem.\! 3.7, to get $(a)$ and
 $(iii)$, mainly \eqref{thm62} and the paragraph before $(ii)$, to get $(b)$).
 Note that if $C$ has seminatural cohomology (i.e. maximal rank and maximal
 corank), then the assumptions of $(a)$ and $(b)$ obviously hold, and we get
 Prop.\! 3.2 of \cite{MDP2}, ch.\! IV, which leads to \cite{MDP2}, ch.\!
 V, Prop.\! 2.1 and to the unobstructedness of $C$ in the case $\diam M
 \leq 2$ (the latter is also proved in \cite{BM}).

 \begin{proposition} \label{gen}
   Let $C$ be any curve in $\proj{3}$ and let $M = \HH_{*}^1(\sI_C)$ and $E =
   \HH_{*}^1(\sO_{C})$. Then
\begin{equation*}
 (a)   {\hspace{1.5 cm}} {_0\!\Hom_R}(I,M) = 0 \ \ {\rm implies} \ \
    \HH_{ \gamma} \cong H(d,g) \ \ {\rm as \ schemes \ at} \  (C)
    {\hspace{1.5 cm}} 
\end{equation*}
\begin{equation*}
(b) {\hspace{1.5 cm}}  {_0\!\Hom_R}(M, E) = 0 \ \ \ {\rm implies} \ \ \
\HH_{\gamma,\rho} \cong \HH_{\gamma}\ {\rm \ as \ schemes \ at} \ (C)\ .
{\hspace{1cm}}  
\ \  
\end{equation*}
\end{proposition}

Finally, we shall in Section 4 see what happens to the unobstructedness of $C$
when we impose on $C$ different conditions of being "general enough". One
result is already now clear, and it points out that the condition $(iii)$ of
Theorem~\ref{munobstr} is the most important one for {\it generic} curves:

\begin{proposition} \label{gen1}
  Let $C$ be a curve in $\proj{3}$, and suppose $C$ is generic in the Hilbert
  scheme $H(d,g)$ and satisfies ${_0\!\Ext_R^2}(M , M) = 0$.  Then $C$ is
  unobstructed if and only if
\begin{equation}
  {_0\!\Hom_R}(I , M) = 0 \ \ {\rm and} \ \ {_0\!\Hom_R}(M , E) = 0 \ .
\end{equation}
\end{proposition}

\begin{proof}  One way is clear from Theorem~\ref{munobstr}. Now suppose $C$
  is unobstructed and generic with postulation $\gamma$ and deficiency $\rho$.
  By generic flatness we see that H$_{\gamma,\rho} \cong \HH_{\gamma} \cong
  \HH(d,g)$ near $C$ from which we deduce an isomorphism of tangent spaces
  T$_{\gamma,\rho} \cong {_0\!\Ext_R^1}(I , I) \cong \HH^0(\sN_C)$. We
  therefore conclude by the exact sequences \eqref{thm62} and \eqref{mseq},
  recalling that $\alpha : {_0\!\Ext_{\mathfrak m}^2}(I , I) \rightarrow
  {_0\!\Ext_R^2}(I , I)$, which appears in \eqref{mseq} for $v = 0$ factorizes
  via ${_0\!\Ext_R^2}(M , M)$, i.e.  $\alpha = 0$.
\end{proof}

\begin{remark} \label{remmunobstr3}
  Combining \eqref{thm62} and \eqref{thm65} we get a {\it surjective\/} map
  $T_{\gamma,\rho} \rightarrow {_0\!\Ext_R^1}(M,M)$ whose kernel
  $T_{\gamma,M}$ is the tangent space of $\HH_{\gamma,M}$ at $(C)$. Now
  dualizing the exact sequence of \eqref{lem25} (for $v = -4$), one proves
  that the surjective map above fits into the exact sequence
\begin{equation}
  k  \rightarrow   {_0\!\Hom_R}(M ,M)  \rightarrow
     {_{-4}\!\Hom_R}(I,E)^{\vee}  \rightarrow   T_{ \gamma , \rho } \rightarrow
     {_0\!\Ext_R^1} (M ,M )   \rightarrow   0   
\end{equation}
and $k \rightarrow {_0\!\Hom_R}(M, M)$ is injective provided $M \neq 0$. We
can use this surjectivity (and some considerations on the obstructions
involved) to give a new proof of the smoothness of the morphism from $\HH_{
  \gamma,\rho}$ to the "scheme" of Rao modules (\cite{MDP1}, Thm. 1.5, page
135). Since ${_{-4}\!\hom}(I,E) = \delta^2(-4)$, cf. \eqref{lem26}, the exact
sequence above also leads to the dimension formula of $\HH_{\gamma,M}$ we
pointed out in Remark~\ref{remeulerI}.
\end{remark}

\section{Sufficient conditions for {\bf obstructedness}}
In this section we will prove that the conditions $(i)$, $(ii)$, $(iii)$ of
Theorem~\ref{munobstr} are both necessary and sufficient for unobstructedness
provided $M$ has diameter one. More generally we are, under the assumption
${_0\!\Ext_R^2}(M,M)=0$ (resp. $\diam M = 1$), able to make explicit conditions
which imply (resp. are equivalent to) obstructedness.  Indeed note that we can
immediately reformulate the first part of Theorem~\ref{munobstr} as
\begin {proposition} \label{propobstr} 
  Let $C$ be a curve in $\proj{3}$, 
  and let ${_0\!\Ext_R^2}(M,M)=0$. If $C$ is obstructed, then (at least) one of
  the following conditions hold  \\[-3mm]

$\ \ \ \ (a) \ \ \ \ \ {_0\!\Hom_R}(I , M) \neq 0 \ \ \ {\rm and} \ \ \
{_0\!\Hom_R}(M , E) \neq 0 \ ,  $ \\[-3mm]

 $\ \ \ \ (b) \ \ \ \ {_{-4}\!\Hom_R}(I , M) \neq 0 \ \ \ {\rm and} \ \ \
  {_0\!\Hom_R}(M , E) \neq 0  \ ,  $  \\[-3mm]

 $\ \ \ \ (c) \ \ \ \  \ {_0\!\Hom_R}(I , M) \neq 0 \ \ \ \ {\rm and}  \ \
{_{-4}\!\Hom_R}(M , E) \neq 0 \ .  $  
\end {proposition}

If $C$ in addition is Buchsbaum, or more generally if the $R$-module $M$
contains ``a Buchsbaum component'', by which we mean that $M$ admits a
decomposition $M=M' \oplus M_{[t]}$ as $R$-{\it modules} where the diameter of
$ M_{[t]}$ is $1$ (i.e. the surjection $M \rightarrow M_{[t]}$ splits as an
$R$-linear map), then we have the following ``converse'' of
Proposition~\ref{propobstr}.
\begin {theorem} \label{obstr} 
  Let $C$ be a curve in $\proj{3}$, let $M = \HH_{*}^1(\sI_C)$ and $E =
  \HH_{*}^1(\sO_{C})$ and suppose ${_0\!\Ext_R^2}(M,M)=0$. Moreover suppose
  there is an $R$-module isomorphism $M \cong M' \oplus M_{[t]}$ where the
  diameter of $ M_{[t]}$ is $1$ and $ M_{[t]}$ supported in degree $t$. Then
  $C$ is obstructed if at least one of the following conditions hold

$ \ \ \ \ (a) \ \ \ \ \  \ {_0\!\Hom_R}(I ,  M_{[t]}) \neq 0 \ \ \ {\rm and}
\ \ \ {_0\!\Hom_R}( M_{[t]} , E) \neq 0 \ , \ \  {\rm or} $  \\[-3mm]

$ \ \ \ \ (b) \ \ \ \ \ {_{-4}\!\Hom_R}(I , M_{[t]}) \neq 0 \ \ \ {\rm and} \
\ \  {_0\!\Hom_R}(M_{[t]},E) \neq 0  \ , \ \  {\rm or}  $   \\[-3mm]

$ \ \ \ \ (c) \ \ \ \ \ \ \ {_0\!\Hom_R}(I ,  M_{[t]}) \neq 0 \ \ \  {\rm and} 
\ {_{-4}\!\Hom_R}( M_{[t]} , E) \neq 0 \ . $
\end{theorem}

Note that if we consider curves obtained by applying Liaison Addition to two
curves where one of them is Buchsbaum of diameter 1, then we always have a
decomposition of $M$ as in Theorem~\ref{obstr} (\cite{MIG}, Thm.\! 3.2.3), see
also \cite{MDP4} for some other cases. Moreover observe that if the module
$L_2$ below has no generators in degree $t$ and $t+4$, then the condition
${_0\!\Ext_R^2}(M,M)=0$ holds if it holds for $M'$, i.e.
${_0\!\Ext_R^2}(M',M')=0$ (Remark~\ref{remmcor}). We get Theorem~\ref{obstr}
immediately from Proposition~\ref{cupprod1} and \ref{cupprod2} which we prove
shortly.

We can state Theorem~\ref{obstr} in terms of the non-triviality of certain
graded Betti numbers of the homogeneous ideal $I = I(C)$. To see this, recall
that once we have a minimal resolution of the Rao module $M$ of free graded
$R$-modules, 
\begin{equation} \label{resoluM}
  0 \rightarrow L_4 \stackrel{\sigma}{\longrightarrow}
    L_3 \rightarrow L_2 \rightarrow L_1
    \rightarrow L_0 \rightarrow M \rightarrow 0 \ \ , \ \
\end{equation}
one may put the unique minimal resolution \eqref{resolu} of the homogeneous
ideal $I$, $0 \rightarrow \oplus_i R(-i)^{\beta_{3,i}} \rightarrow \oplus_i
R(-i)^{\beta_{2,i}} \rightarrow \oplus_i R(-i)^{\beta_{1,i}} \rightarrow I
\rightarrow 0 $, in the following form
\begin{equation} \label{resoluMI}
  0 \rightarrow L_4 \stackrel{\sigma \oplus 0}{\longrightarrow} L_3 \oplus F_2
    \rightarrow F_1 \rightarrow I \rightarrow 0 \ , \ 
\end{equation}
i.e. where the composition of $L_4 \rightarrow L_3 \oplus F_2$ and the natural
projection $L_3 \oplus F_2 \rightarrow F_2 $ is zero (\cite{R}, Theorem 2.5).
Note that any minimal resolution of $I$ of the form \eqref{resoluMI} has
well-defined modules $F_2$ and $F_1$. In particular $ F_1= \oplus_i
R(-i)^{\beta_{1,i}}$. Moreover applying $\Hom(-,M)$ to \eqref{resoluM} we get
a minimal resolution of ${\Ext_R^4}(M,R) \cong {\Ext_R^4}(M' \oplus M_{[t]},
R) \cong {\Ext_R^4}(M', R) \oplus M_{[t]}(2t+4)$ from which we see that $L_4$
contains $R(-t-4)^r$ as a direct summand where $r = \dim_k M_{[t]}$. Put
\begin{equation} \label{bettiMI}
 L_4 \cong L_4' \oplus R(-t-4)^r \ , \ F_2 \cong P_2 \oplus R(-t-4)^{b _1}
    \oplus R(-t)^{b _2} \ , \ F_1 \cong  P_1 \oplus 
    R(-t-4)^{a_1} \oplus R(-t)^{a_2}
\end{equation}
where $P_i$, for $i = 1,2$ are supposed to contain no direct factor of degree
$t$ and $t+4$. So $a_1$ and $a_2$ are exactly the first graded Betti number of
$I$ in the degree $t+4$ and $t$ respectively, while $b_1$ and $r$ (resp.
$b_2$) are less than or equal to the corresponding Betti number of $I$ in
degree $t+4$ (resp. $t$) because $L_4'$ and $L_3$ might contribute to the
graded Betti numbers. If, however, $M$ is of diameter $1$ (and $M \cong
M_{[t]}$), then $L_4'=0$ and the generators of $L_3$ sit in degree $t+3$. In
this case $b_i$ and $r$ are exactly equal to the corresponding graded Betti
numbers in the minimal resolution \eqref{resolu}. Now Theorem~\ref{obstr}
translates to

\begin {corollary} \label{mcor}
  Let $C$ be a curve in $\proj{3}$, let ${_0\!\Ext_R^2}(M,M)=0$ and  suppose $M
  \cong M' \oplus M_{[t]}$ as $R$-{\it modules} where the diameter of $
  M_{[t]}$ is $1$ and supported in degree $t$. Then $C$ is obstructed if 
  $$
  a_2 \cdot b_1 \neq 0 \ \ \ {\rm or} \ \ \ a_1 \cdot b_1 \neq 0 \ \ \ {\rm
    or} \ \ \ a_2 \cdot b_2 \neq 0 \ . $$
\end {corollary}

This leads to one of the main Theorems of this paper, which solves the problem
of characterizing obstructedness in the diameter 1 case (raised in \cite{EF1})
completely.

\begin{theorem} \label{mainthm}
  Let $C$ be a curve in $\proj{3}$ whose Rao module $M \neq 0$ is of diameter
  $1$ and concentrated in degree $c$, and let $\beta_{1,c+4}$ and
  $\beta_{1,c}$ (resp.  $\beta_{2,c+4}$ and $\beta_{2,c}$) be the number of
  minimal generators (resp. minimal relations) of $I$ of degree $c+4$ and $c$
  respectively. Then $C$ is obstructed if and only if
  $$
  \beta_{1,c} \cdot \beta_{2,c+4} \neq 0 \ \ \ {\rm or} \ \ \ \beta_{1,c+4}
  \cdot \beta_{2,c+4} \neq 0 \ \ \ {\rm or} \ \ \ \beta_{1,c} \cdot
  \beta_{2,c}\neq 0 \ . $$
  Moreover if $C$ is {\bf unobstructed} and $M$ is
  $r$-dimensional (i.e. $r=\beta_{3,c+4}$), then the dimension of the Hilbert
  scheme $\HH(d,g)$ at $(C)$ is
\begin{equation*}
 \dim_{(C)}\HH(d,g) = 4d + \delta^2(0) + r(\beta_{1,c+4} + \beta_{2,c}).
\end{equation*}
\end{theorem}

\begin{proof}[Proof (of Corollary~\ref{mcor})] \ In the sequel we frequently
  use the triviality of the module structure of $M_{[t]}$ ($\mathfrak m \cdot
  M_{[t]} = 0$). Now applying ${_v\!\Hom_R}(- ,M_{[t]})$ to the minimal
  resolution \eqref{resoluMI} we have by \eqref{bettiMI},
\begin{equation} \label{cor61}
\ {_0\!\hom_R}(I,M_{[t]}) = ra_2 \ \  {\rm and}  \ \
  {_{-4}\!\hom_R}(I, M_{[t]}) = ra_1 \ .
\end{equation}
Moreover note that the assumption ${_0\!\Ext_R^2}(M,M)=0$ implies
${_{-4}\!\Ext_R^2}(M,M)=0$ by \eqref{thm67} and hence
${_v\!\Ext_R^2}(M_{[t]},M)=0$ for $v=0$ and $-4$ by the split $R$-linear map
$M \rightarrow M_{[t]}$. By the duality \eqref{mduality} and the spectral
sequence \eqref{mspect} (which converges to ${_v\!\Ext_{\mathfrak
    m}^3}(M_{[t]},I))$ we therefore get an exact sequence
\begin{equation} \label{cor62}  
0 \rightarrow {_v\!\Ext_R^1}(M_{[t]},M) \rightarrow
        {_{-v-4}\!\Ext_R^1}(I,M_{[t]})^{\vee} \rightarrow
        {_v\!\Hom_R}(M_{[t]}, E) \rightarrow 0
\end{equation}
for $v=0$ and $-4$. Since $ {_{v}\!\Ext_R^1}(M_{[t]}, M)^{\vee} \cong
{_{-v-4}\!\Ext_R^3}(M,M_{[t]})$ by \eqref{mduality} and \eqref{mspect} and
since we have ${_{-v-4}\!\Ext_R^3}(M,M_{[t]}) \cong
{_{-v-4}\!\Hom_R}(L_3,M_{[t]})$ by \eqref{resoluM}, we get $
{_{v}\!\Ext_R^1}(M_{[t]}, M) \cong {_{-v-4}\!\Hom_R}(L_3,M_{[t]})^{\vee}$.
Interpreting ${_{-v-4}\!\Ext_R^1}(I ,M_{[t]})$ similarly via the minimal
resolution \eqref{resoluMI} of $I$, we get $ {_v\!\Hom_R}(M_{[t]}, E) \cong
{_{-v-4}\!\Hom_R}(F_2,M_{[t]})^{\vee}$ for $v=0$ and $-4$ and hence
\begin{equation} \label{cor63}
  {_0\!\hom_R}(M_{[t]}, E) = rb_1 \quad {\rm and} \quad
  {_{-4}\!\hom_R}(M_{[t]}, E) = rb_2
\end{equation}
by \eqref{bettiMI} and we conclude easily since $r \neq 0$.
\end{proof}

\begin{remark} \label{remmcor}
  For later use, note that ${_{v}\!\Ext_R^2}(M_{[t]}, M))^{\vee} \cong
  {_{-v-4}\!\Ext_R^2}(M,M_{[t]}) \cong {_{-v-4}\!\Hom_R}(L_2,M_{[t]})$. Hence
  if we assume the latter group to vanish (instead of assuming
  ${_0\!\Ext_R^2}(M,M)=0$), we get \eqref{cor62} and \eqref{cor63} for this
  $v$. In particular if $ {_{v}\!\Ext_R^2}(M_{[t]},M)=0$ for $v=0$ and $-4$,
  then \eqref{cor62} and \eqref{cor63} hold, as well as $ {_0\!\Ext_R^2}(M,M)
  \cong {_0\!\Ext_R^2}(M',M')$ because ${_0\!\Ext_R^2}(M', M_{[t]}) \cong
  {_{-4}\!\Ext_R^2}( M_{[t]},M')^{\vee}=0$.
\end{remark}

\begin{proof}[Proof (of Theorem~\ref{mainthm})] \ Combining
  Proposition~\ref{propobstr} and Corollary~\ref{mcor} we immediately get the
  first part of the Theorem.  Moreover since we by Remark~\ref{remunobstr}
  have
\begin{equation*} 
 \dim_{(C)}\HH(d,g) = 4d + \delta^2(0) + {_{-4}\!\hom_R}(I,M) +
  {_{-4}\!\hom_R}(M , E)\ , 
\end{equation*}
we conclude by \eqref{cor61} and \eqref{cor63}.
\end{proof}

To prove Theorem~\ref{obstr} the following key proposition is useful. As
Fløystad points out in \cite{F}, if the image of the cup product
$<\lambda,\,\lambda>\ \in \Ext_{\sO_{\proj{}}}^2(I_C,I_C)$, $\lambda \ \in
\Ext_{\sO_{\proj{}}}^1(I_C,I_C)$, maps to a non-zero element $\bar o \in
{_0\!\Hom_R}(I,E)$ via the right vertical map of \eqref{cupsquare} below, then
$C$ is obstructed. He makes several nice contributions to calculate $\bar o$,
especially when $M$ is a complete intersection (e.g. \cite{F}, Prop.\! 2.13,
from which Proposition~\ref{cupprod1} is an easy consequence, and \cite{F}, $
§ 5$), see also \cite{MDP3}, $ § 3$ for further calculations
and Laudal (\cite{L2}, $ § 2$) for the theory of cup and Massey products. In
general it is, however, quite difficult to prove that $\bar o \neq 0$, while
the non-vanishing of the natural composition
$$
{_0\!\Hom_R}(I,M) \times {_0\!\Hom_R}(M , E) \longrightarrow {_0\!\Hom_R}(I ,
E)$$ is easier to handle. This is the benefit of
Proposition~\ref{cupprod1} which we prove by using that $\alpha = 0$ in
\eqref{mseq}.

\begin {proposition} \label{cupprod1}  
  Let $C$ be a curve in $\proj{3}$, let $M = \HH_{*}^1(\sI_C)$ and $E =
  \HH_{*}^1(\sO_C)$ and suppose ${_0\!\Ext_R^2}(M , M) = 0$. If the obvious
  morphism
\begin{equation*}
  {_0\!\Hom_R}(I,M) \times {_0\!\Hom_R}(M , E) \longrightarrow  {_0\!\Hom_R}(I ,
  E) 
\end{equation*}
(given by the composition) is non-zero, then $C$ is obstructed. In particular
if $M$ admits a decomposition $M=M' \oplus M_{[t]}$ as $R$-{\it modules} where
the diameter of $ M_{[t]}$ is $1$, then $C$ is obstructed provided
$$  
 {_0\!\Hom_R}(I, M_{[t]}) \neq 0 \ \ {\rm and} \ \ \ {_0\!\Hom_R}( M_{[t]}, E)
 \neq 0 $$ 
\end{proposition}

\begin{proof} It is well known (cf. \cite{L1}) that if the Yoneda pairing
(inducing the cup product)
\begin{equation*}
  <-,-> \ \ : \ \ \ \Ext_{\sO_{\proj{}}}^1(\sI_C , \sI_C) \times
    \Ext_{\sO_{\proj{}}}^1(\sI_C, \sI_C) \rightarrow
    \Ext_{\sO_{\proj{}}}^2(\sI_C,\sI_C), 
\end{equation*}
given by composition of resolving complexes, satisfies $<\lambda,\lambda> \neq
0$ for some $\lambda$, then $C$ is obstructed. If we let $p_1 :
\Ext_{\sO_{\proj{}}}^1(\sI_C,\sI_C) \rightarrow {_0\!\Hom_R}(I, M)$ and $p_2:
\Ext_{\sO_{\proj{}}}^1(\sI_C ,\sI_C) \rightarrow {_0\!\Hom_R}(M, E)$ be the maps
induced by sending an extension onto the corresponding connecting
homomorphisms, then $<$-,-$>$ fits into a commutative diagram

\begin{equation} \label{cupsquare}
\begin{array}{ccccccccc}
    &  \Ext_{\sO_{\proj{}}}^1(\sI_C,\sI_C) & \times &
      \Ext_{\sO_{\proj{}}}^1(\sI_C,\sI_C) 
      & \longrightarrow &  \Ext_{\sO_{\proj{}}}^2(\sI_C,\sI_C) &  \\
    &  \ \downarrow p_1 \ & & \ \ \downarrow p_2 \ & & \downarrow &  &  \\
   & _0\!\Hom_R(I,M) & \times & _0\!\Hom_R(M,E) & \longrightarrow &
    _0\!\Hom_R(I,E)   &  &
\end{array}
\end{equation} 
where the lower horizontal map is given as in Proposition~\ref{cupprod1}. By
\eqref{mseq}, ${_0\!\Ext_R^1}(I , I)=\ker p_1$, and $p_1$ is surjective
because $\alpha = 0$ for $v = 0$. Moreover since the composition
${_0\!\Ext_R^1}(I , I) \hookrightarrow \Ext^1(\sI_C,\sI_C) \rightarrow
{_0\!\Hom_R}(M , E)$ is surjective by the important sequence \eqref{thm62},
there exists $(\lambda_1 ,\lambda_2) \in \Ext^1(\sI_C,\sI_C) \times
{_0\!\Ext_R^1}(I , I)$ such that the composed map
$p_2(\lambda_2)p_1(\lambda_1)$ is non-zero by assumption. Using $\lambda_2 \in
{_0\!\Ext_R^1}(I , I) = \ker p_1$, we get
\begin{equation*}
 p_2(\lambda_1 + \lambda_2)p_1(\lambda_1 + \lambda_2) =
    p_2(\lambda_1)p_1(\lambda_1) + p_2(\lambda_2)p_1(\lambda_1)
\end{equation*}
i.e. either $<\lambda_1 + \lambda_2\ ,\ \lambda_1 + \lambda_2>$ or $<\lambda_1
\ ,\ \lambda_1>$ are non-zero, and $C$ is obstructed.

Finally suppose the two last mentioned $\Hom$-groups of
Proposition~\ref{cupprod1} are non-vanishing.  Then there exists a map $\psi
\in {_0\!\Hom_R}(M_{[t]}, E)$ such that $\psi(m) \neq 0$ for some $m \in
(M_{[t]})_t$. Since $M_{[t]}$ has diameter $1$, we get
${_0\!\Hom_R}(I,M_{[t]}) \cong {_0\!\Hom_R}(R(-t)^{a_2},M_{[t]}) \cong
(M_{[t]})_t^{a_2}$ by \eqref{resoluMI} and \eqref{bettiMI}, and we have $a_2 >
0$. Hence there is a map $\phi' \in {_0\!\Hom_R}(I,M_{[t]})$ such that
$\phi'(1,0,...,0) = m$ where $(1,0,...,0)$ an $a_2$-tuple. Since $
{_0\!\Hom}(I, M) \rightarrow {_0\!\Hom}(I, M_{[t]}) $ is surjective by the
existence of the $R$-split morphism $p:M \rightarrow M_{[t]}$ there is an
element $\phi \in {_0\!\Hom}(I, M)$ which maps to $\phi'$. Since the
composition $\psi \phi' = \psi p \phi$
maps to a non-trivial element of ${_0\!\Hom_R}(I , E)$ by construction, we
conclude by the first part of the proof.
\end{proof}

\begin{remark}  \label{remcupprod1} 
  Let $C$ be a curve in $\proj{3}$ whose Rao module has {\it diameter\/} $1$.
  From \eqref{mseq} and \eqref{thm62}, cf.  the proof above, we see at once
  that ${_0\!\Hom_R}(I , M) \neq 0$ and ${_0\!\Hom_R}(M , E) \neq 0$ if and
  only if we have the following strict inclusions of tangent spaces
\begin{equation} \label{incl}
 T_{\gamma,\rho} \varsubsetneq {_0\!\Ext_R^1}(I,I) \varsubsetneq \HH^0(\sN_C)
\end{equation}
where ${_0\!\Ext_R^1}(I , I)$ is the tangent space of the Hilbert scheme of
constant postulation $\HH_{\gamma}$ at $(C)$. By Proposition~\ref{cupprod1},
{\it $C$ is obstructed if \eqref{incl} holds}. If $M \cong k$, this conclusion
follows also from \cite{MDP1}, ch.\! $ X$, Prop.\! 5.9, or from  \cite{MDP3}.
\end{remark}

Along the same lines we are able to generalize a result of Walter \cite{W1}. 
If the diameter of $M$ is $1$ and if ${_0\!\Hom_R}(I , M) = 0$, then Walter
proves Proposition~\ref{cupprod2}$(a)$ below and he computes the completion of
$\sO_{\HH(d,g),(C)}$ in detail.  The first part of Proposition~\ref{cupprod1}
and \ref{cupprod2}, however, requires only ${_0\!\Ext_R^2}(M , M) = 0$. This
vanishing condition, which one may show is invariant under linkage (by using
\eqref{prop41}), holds for instance if the diameter of $M$ is less or equal
$2$, or if $M$ is generic of diameter $3$ and the scheme of Rao modules is
irreducible (cf. \cite{MDP2}).

\begin{proposition}  \label{cupprod2}
  Let $C$ be a curve in $\proj{3}$, let $M = \HH_{*}^1(\sI_C)$, $E =
  \HH_{*}^1(\sO_C)$ and let ${_0\!\Ext_R^2}(M , M) =  0$. 

 (a) If the obvious morphism
\begin{equation*}
  {_{-4}\!\Hom_R}(I,M) \times {_0\!\Hom_R}(M , E) \longrightarrow
  {_{-4}\!\Hom_R}(I , E)
\end{equation*}
(given by the composition) is non-zero, then $C$ is obstructed. In particular
if $M$ admits a decomposition $M=M' \oplus M_{[t]}$ as $R$-{\it modules} where
the diameter of $ M_{[t]}$ is $1$, then $C$ is obstructed provided
$$
{_{-4}\!\Hom_R}(I,M_{[t]}) \neq 0 \ \ {\rm and} \ \ \ {_0\!\Hom_R}(M_{[t]},
E) \neq 0 $$

 (b) If the morphism
\begin{equation*}
  {_0\!\Hom_R}(I,M) \times {_{-4}\!\Hom_R}(M , E) \longrightarrow
  {_{-4}\!\Hom_R}(I,E)  
\end{equation*}
(given by the composition) is non-zero, then $C$ is obstructed. In particular
if $M$ admits a decomposition $M=M' \oplus M_{[t]}$ as $R$-{\it modules} where
the diameter of $ M_{[t]}$ is $1$, then $C$ is obstructed provided
$$
{_0\!\Hom_R}(I,M_{[t]}) \neq 0 \ \ {\rm and} \ \ \ {_{-4}\!\Hom_R}(M_{[t]},
E) \neq 0 $$
\end{proposition}

\begin{proof} {\it Step} $1$.  In Step 1 we give a full proof of $(a)$ under
  the extra temporary assumption $M_{-4} = 0$. Denote by $p_2'$ the
  restriction of $p_2$ (see \eqref{cupsquare}) to ${_0\!\Ext_R^1}(I,I)$ via
  the natural inclusion ${_0\!\Ext_R^1}(I,I) \hookrightarrow
  \Ext^1(\sI_C,\sI_C)$ and consider the commutative diagram
\begin{equation} \label{cupsquare2}
\begin{array}{ccccccccc}
 <-,->_0 \ : &  {_0\!\Ext_R^1}(I,I) & \times &  {_0\!\Ext_R^1}(I,I) 
      & \longrightarrow &  {_0\!\Ext_R^2}(I,I) &  \\
    &  \ \uparrow  \ & & \ \ \downarrow p_2 ' \ & & \downarrow i &  &  \\
   &  T_{\gamma,\rho} & \times & _0\!\Hom_R(M,E) & \longrightarrow &
    {_0\!\Ext_R^1}(M,E)   &  &
\end{array}
\end{equation} 
where $<-,->_0$ is the Yoneda pairing. Indeed the restriction of
${_0\!\Ext_R^1}(I , I)$ to the subspace $T_{\gamma,\rho}$ in
\eqref{cupsquare2} makes the lower horizontal arrow well-defined in the
commutative diagram above because of the natural map
$T_{\gamma,\rho}\rightarrow {_0\!\Ext_R^1}(M , M)$ of
Remark~\ref{remmunobstr3}. Due to the exact sequence \eqref{thm62}, continued
as in Remark~\ref{remmunobstr}, the map $p_2'$ is surjective and $i$ is
injective by the assumption ${_0\!\Ext_R^2}(M, M) = 0$. Hence the pairing
$<-,->_0$ factorizes via
\begin{equation} \label{prop82}
 \varphi'\ : \ T_{\gamma,\rho}  \times  _0\!\Hom_R(M,E)  \longrightarrow 
     {_0\!\Ext_R^2}(I,I)   
\end{equation}
and vanishes if we restrict $\varphi'$ to ${_{-4}\!\Hom_R}(I, E)^{\vee}
\times {_0\!\Hom_R}(M , E)$ via the map of Remark~\ref{remmunobstr3}. (using
the identity on ${_0\!\Hom_R}(M , E)$, because ${_{-4}\!\Hom_R}(I , E)^{\vee}$
maps to zero in $ {_0\!\Ext_R^1}(M , M)$. 

To prove $(a)$ it suffices to prove $<\lambda\, ,\, \lambda >_0 \, \neq 0$ for
some $\lambda$. We do this, we {\it claim that \/}there is another pairing
$\varphi \neq 0$, which commutes with $<-,->_0$, and which essentially
corresponds to the restriction of $\varphi'$ above except for the exchange of
variables, i.e.
\begin{equation} \label{prop83}
 \varphi\ : \ _0\!\Hom_R(M,E) \times {_{-4}\!\Hom_R}(I, E)^{\vee}
     \longrightarrow  {_0\!\Ext_R^2}(I,I)   
\end{equation}
(Since $T_{\gamma,M} = \coker({_0\!\Hom_R}(M, M) \rightarrow
{_{-4}\!\Hom_R}(I, E)^{\vee})$ by Remark~\ref{remmunobstr3}, we can continue
the arguments below to see that the map $\varphi$ of \eqref{prop83} extends to
a somewhat more naturally defined pairing ${_0\!\Hom_R}(M , E) \times
T_{\gamma,M} \rightarrow {_0\!\Ext_R^2}(I , I)$, but this observation does not
really effect the proof). Now, to prove the claim there is, as in
\eqref{cupsquare}, a commutative diagram
\begin{equation*} 
\begin{array}{ccccccccc}
    &  _{-4}\!\Ext_{\mathfrak m}^2(I,I) & \times &
      _0\!\Ext_R^1(I,I) 
      & \longrightarrow &  _{-4}\!\Ext_{\mathfrak m}^3(I,I) &  \\
    &  \ \downarrow \cong \ & & \ \ \downarrow p_2' \ & & \downarrow &  &  \\
   & _{-4}\!\Hom_R(I,M) & \times & _0\!\Hom_R(M,E) & \longrightarrow &
    _{-4}\!\Hom_R(I,E)   &  & \\
\end{array}
\end{equation*} 
where two of the vertical arrows are given by the spectral sequence
\eqref{mspect} (cf. \eqref{lem25}) and where the lower pairing is the {\it
  non-vanishing\/} map of Proposition~\ref{cupprod2}. Dualizing, we get the
commutative diagram
\begin{equation*} 
\begin{array}{ccccccccc}
 &  {_0\!\Ext_R^1}(I,I) & \times &  {_{-4}\!\Ext_{\mathfrak
      m}^3}(I,I)^{\vee}  
      & \longrightarrow &  {_{-4}\!\Ext_{\mathfrak m}^2}(I,I)^{\vee} &  \\
    &  \ \downarrow p_2'  \ & & \ \ \uparrow  \ & & \uparrow \cong &  &  \\
   &   _0\!\Hom_R(M,E) & \times &  {_{-4}\!\Hom_R}(I,E)^{\vee} &
      \longrightarrow &  {_{-4}\!\Hom_R}(I,M)^{\vee}   &  &
\end{array}
\end{equation*} 
where the {\it non-vanishing} lower arrow can be identified with the map
$\varphi$ of \eqref{prop83}. Using the duality \eqref{mduality}, we see that
$\varphi$ commutes with the Yoneda pairing $<-,->_0$, and the claim follows
easily.

Now since $\varphi \neq 0$ and $p_2'$ is surjective, there exists $(\lambda_2,
\lambda_1) \in {_0\!\Hom_R}(M , E) \times {_{-4}\!\Hom_R}(I,E)^{\vee}$ and
$\lambda_2'\in {_0\!\Ext_R^1}(I,I)$ such that $p_2'(\lambda_2')\ =\ \lambda_2$
and such that $<\lambda_2' ,\lambda_1 >_0 \ =\ \varphi(\lambda_2,\lambda_1)
\neq 0$. Note that $<\lambda_1 ,\lambda>_0\ =\ 0$ for any $\lambda \in
{_0\!\Ext_R^1}(I , I)$ because $<\lambda_1 ,\lambda>_0\ = \ \varphi'(\lambda_1
,p_2'(\lambda)) = 0$ by \eqref{prop82}. It follows that
\begin{equation*}
  <\lambda_1 + \lambda_2',\lambda_1
    + \lambda_2'>_0 \ =\ <\lambda_2',\lambda_1>_0 + <\lambda_2',\lambda_2'>_0
\end{equation*}
i.e. either $<\lambda_1 + \lambda_2',\lambda_1 + \lambda_2'>_0$ or
$<\lambda_2',\lambda_2'>_0$ are non-zero. Finally since the map $\alpha$ of
\eqref{mseq} factors via ${_0\!\Ext}_R^2(M ,M)$ for $v = 0$, it follows that
the map ${_0\!\Ext_R^2}(I , I) \rightarrow Ext^2(\sI_C,\sI_C)$ is injective and
maps obstructions to obstructions, i.e. the Yoneda pairing $<-,->_0$ and the
corresponding pairing $<-,->$ of \eqref{cupsquare} commute and vanish
simultaneously. $C$ is therefore obstructed.

{\it Step\/} 2. To prove $(b)$ we use Step 1 and Proposition~\ref{mlink}.
Indeed let $C$ be a curve as in $(b)$ and let $Y \supseteq C$ be a complete
intersection of two surfaces of degrees $f$ and $g$ such that the {\it
  conditions\/} of Proposition~\ref{mlink} hold and such that
$\HH^1(\sI_C(f+g)) = 0$, $\HH^1(\sO_{C}(f-4)) = 0$ and $\HH^1(\sO_C(g-4)) = 0$
(such $Y$ exists). Then we claim that the corresponding linked curve $C'$
satisfies the conditions given in Step 1. Indeed slightly extending
Remark~\ref{remmunobstr}, we have
\begin{equation*}
  {_0\!\Hom_R}(I(C), M(C)) \cong {_0\!\Hom_R}(M(C'), E(C'))
\end{equation*}
\begin{equation} \label{prop84}
  {_{-4}\!\Hom_R}(M(C) , E(C)) \cong {_{-4}\!\Hom_R}(I(C'), M(C'))
\end{equation}
\begin{equation*}
  {_{-4}\!\Hom_R}(I(C ), E(C)) \cong {_{-4}\!\Hom_R}(I(C)/I(Y), E(C)) \cong
        {_{-4}\!\Hom_R}(I(C')/I(Y), E(C'))
\end{equation*}
and we get the claim because ${_{-4}\!\Hom_R}(I(C')/I(Y) , E(C')) \rightarrow
{_{-4}\Hom_R}(I(C') , E(C'))$ is injective and $\HH^1(\sI_C(f+g)) \cong
\HH^1(\sI_{C'}(-4))$. It follows that $C'$ is obstructed by Step 1, and so is
$C$ by Proposition~\ref{mlink}. Moreover if $M=M' \oplus M_{[t]}$ and the
diameter of $M_{[t]}$ is 1, we conclude easily by arguing as in the very end
of the proof of Proposition~\ref{cupprod1}.

{\it Step\/} 3.  Finally using the same idea as in Step 2, we prove that $(b)$
and Proposition~\ref{mlink} imply $(a)$. Indeed by Proposition~\ref{mlink} we
can see that $(a)$ and $(b)$ are equivalent by making a suitable linkage, and
the proof is complete.
\end{proof}

Focusing on the Hilbert scheme with constant postulation, $ \HH_{\gamma}$, we
have the following result, quite similar to Theorem~\ref{mainthm}.

\begin{proposition} \label{remmainthm}
  Let $C$ be a curve in $\proj{3}$ whose Rao module $M \neq 0$ is of diameter
  $1$ and concentrated in degree $c$, and let $\beta_{1,c+4}$ and
  $\beta_{1,c}$ (resp.  $\beta_{2,c+4}$ and $\beta_{2,c}$) be the number of
  minimal generators (resp. minimal relations) of degree $c+4$ and $c$
  respectively.  Suppose also $M_{-4} = 0$. Then
  $$\HH_{\gamma}\ is \ singular \ at \ (C) \ if \ and \ only \ if \ 
  \beta_{1,c+4} \cdot \beta_{2,c+4}\neq 0 \ . $$
  Moreover if $\ \HH_{\gamma}$
  is smooth at $(C)$ and $M$ is $r$-dimensional (i.e. $r=\beta_{3,c+4}$), then
\begin{equation*}
 \dim_{(C)}\HH_{\gamma} = 4d + \delta^2(0) + r(\beta_{1,c+4} + \beta_{2,c}-
 \beta_{1,c}) \ .
\end{equation*}  
\end {proposition}

 \begin{proof} 
   Since the tangent space, resp.  the obstructions, of $\HH_{\gamma}$ at $C$
   is ${_0\!\Ext_R^1}(I, I)$, resp.  sit in ${_0\!\Ext_R^2}(I,I)$, cf. the
   proof of $(i)$ in Theorem~\ref{munobstr}, we have by Step 1 of the proof
   above that $\HH_{\gamma}$ is not smooth at $(C)$ provided $M_{-4} = 0$ and
   the conditions of Proposition~\ref{cupprod2}$(a)$ hold. Hence if $
   \beta_{1,c+4} \cdot \beta_{2,c+4}\neq 0 \ $, it follows from \eqref{cor61}
   and \eqref{cor63} that $\HH_{\gamma}$ is singular at $(C)$ (only for this
   way here we need the assumption $M_{-4} = 0$).  For the converse, suppose
   $\beta_{1,c+4}=0$.  Then ${_0\!\Ext_R^2}(I,I) \cong {_{-4}\!\hom_R}(I,
   M)^{\vee}=0 $ and if $\beta_{2,c+4}=0$, we get by \eqref{cor63} and
   Proposition~\ref{gen} an isomorphism between $\HH_{\gamma, \rho}$ and
   $\HH_{\gamma}$ at $(C)$. The former scheme is smooth because
   ${_0\!\Ext_R^2}(M , M) = 0$, and we get the smoothness of the latter.
   Finally to see the dimension we use $\chi(\sN_C)=4d$ and \eqref{mseq} with
   $\alpha = 0$ for $v=0$ to get $$
   {_0\!\ext_R^1}(I,I) = 4d + h^1(\sN_C)- {_0\!\hom_R}(I,M)\ , $$
   and we conclude by \eqref{thm51}, \eqref{cor61} and
   \eqref{cor63}.
\end{proof}

\begin{remark} \label{genmcor} 
  Corollary~\ref{mcor} admits the following generalization. Instead of
  assuming the $R$-module isomorphism $M \cong M' \oplus M_{[t]}$, we suppose
  that $M$ contains a minimal generator $T$ of degree $t$ and we replace $a_i
  \neq 0$ by the surjectivity of a certain non-trivial map as follows. Let $M
  \twoheadrightarrow M \otimes_R k \twoheadrightarrow k(-t)$ and $\eta(T) :
  {_v\!\Hom_R}(I , M) \rightarrow {_v\!\Hom_R}(I , k(-t))$ be maps induced by
  $T$. Note that $M \twoheadrightarrow M \otimes_R k $ is not necessarily a
  split $R$-homomorphism. So if $F_{t-v}$ is a minimal generator of $I$ of
  degree $t-v$ (inducing maps $R(-t+v) \hookrightarrow I$ and $\tau(F_{t-v}):
  \ {_v\!\Hom_R}(I,k(-t)) \rightarrow {_v\!\Hom_R}(R(-t+v) ,k(-t)) \cong k$),
  we just suppose the surjectivity of the composition $\tau(F_{t-v}) \eta(T)$
  for $v= 0$ (resp. $-4$) instead of $a_2 \neq 0$ (resp.  $a_1 \neq 0$), to
  get a generalization of Corollary~\ref{mcor}. Hence if
 \begin{equation*}
 \tau(F_{t}) \eta(T)\ {\rm \ is \
 surjective \ for \ some \ minimal \ generator\ } F_{t} \ {\rm \ of \ } I, 
 {\rm \  and\ } b_1 \neq 0 {\rm \ or } \  b_2 \neq  0  \ , \ \ OR 
\end{equation*}
\begin{equation*}
 {\rm \ if \ }  \tau(F_{t+4}) \eta(T)  \ {\rm \ is \  surjective \ for \ some
 \ minimal \ 
 generator\ } F_{t+4} \ {\rm \ of \ }  I, \ {\rm \ and} \  b_1 \neq 0 \ , 
\end{equation*}
then $C$ is obstructed. There is no real change in the proof. Indeed looking
to the very final part of Proposition~\ref{cupprod1} and to the proof of
Corollary~\ref{mcor}, noting that we don't need the surjectivity of $
{_{-v-4}\!\Ext_R^1}(I,k(-t))^{\vee} \rightarrow {_v\!\Hom_R}(k(-t), E) $ in
\eqref{cor62} (where we have replaced $M_{[t]}$ by $k(-t)$), we get the
result. Finally note that it is easy to see that $\tau(F_{t-v}) \eta(T)$ is
surjective if the row in the matrix of relations (i.e. the middle arrow) of
\eqref{resolu} which corresponds to $F_{t-v}$, maps $M_t$ to zero. If $$\min
\{i > t-v \ \arrowvert \ \beta_{2,i} \neq 0 \} > c-t \ , $$
then the entries
of this row map $M_t$ onto $M_{c+j}$ for $j > 0$, i.e onto zero, and we have
the mentioned surjectivity. This surjectivity holds in particular if $t= c$
(and $L_4$ contains generators of degree $c+4$, as always).
\end{remark}

\begin{remark} \label{3Yoneda} 
  We have by Proposition~\ref{cupprod1} and \ref{cupprod2} the following three
  Yoneda pairings
 \begin{equation*}
  {_0\!\Hom_R}(I , M) \times {_0\!\Hom_R}(M, E)\longrightarrow
  {_0\!\Hom_R}(I,E) 
\end{equation*}
\begin{equation*}
  {_0\!\Hom_R}(I,M) \times {_{-4}\!\Hom_R}(I, E)^{\vee}  \longrightarrow
  {_{-4}\!\Hom_R}(M, E)^{\vee}  
\end{equation*}
\begin{equation*}
  {_0\!\Hom_R}(M, E) \times {_{-4}\!\Hom_R}(I , E)^{\vee}  \longrightarrow
  {_{-4}\!\Hom_R}(I , M)^{\vee}  
\end{equation*}
To illustrate, in a diagram, how the right hand sides contribute to
$\HH^1(\sN_C)$, we suppose ${_0\!\Ext_R^i}(M , M) = 0$ for $i \geq 2$ to
simplify. Then recall that ${_0\!\Ext_R^2}(I, M) = 0$ and ${_0\!\Ext_R^1}(I,M)
\cong {_{-4}\!\Hom}(M,E)^{\vee}$ by Remark~\ref{remunobstr}. Now \eqref{mseq}
(resp. \eqref{lem25}) leads to the exactness of the horizontal (resp.
vertical, with injective upper downarrow and surjective lower downarrow)
sequence in the diagram
\begin{equation*} \label{diagramN}
\begin{array}{cccccl}
  &&&  _0\!\Ext_R^1(I,M) & \cong & {_{-4}\!\Hom}(M,E)^{\vee}   \\[1mm]
  &&& \downarrow \\[1mm]
  0 \ \ \longrightarrow & {_0\!\Ext_R^2}(I,I) &  \longrightarrow \ \ \ \
  \HH^1(\sN_C) \hspace{0.3 cm} \ \longrightarrow & {_0\!\Ext_{\mathfrak
      m}^3}(I,I) \hspace{0.1 cm} &  \longrightarrow & {\hspace{0.1 cm}} 0
  \\[1mm] 
  & \downarrow \cong  \ & & \downarrow   \\[1mm]
  & _{-4}\!\Hom_R(I,M)^{\vee}  & &    _0\!\Hom_R(I,E)
\end{array}
\end{equation*} 
\end{remark}

We will end this section by showing that there exists smooth connected space
curves in any of the three cases $(a)$, $(b)$ and $(c)$ of
Theorem~\ref{obstr}. The case $(b)$ is treated in \cite{W1}, where Walter
manages to find obstructed curves of maximal rank (see also \cite{BKM}). These
curves make $\HH_{\gamma}$ singular as well (Proposition~\ref{remmainthm}). By
linkage we can transfer the result in \cite{W1} to the case $(c)$ and we get
the existence of obstructed curves of maximal corank, whose local ring
$\sO_{\HH(d,g),(C)}$ can be described exactly as in \cite{W1}. However, since
we in the next section will see that a sufficiently general curve of
$\HH_{\gamma, \rho}$ does not verify neither $(b)$ nor $(c)$, the case $(a)$
deserves special attention. We shall now see that there exist many smooth
connected curves satisfying the conditions $(a)$.

\begin{example} \label{ex1} We claim that for any triple $(r,a_2,b_1)$ of
  positive integers there exists a smooth connected curve $C$ with minimal
  resolution as in \eqref{resoluMI} and \eqref{bettiMI} and $\diam M(C) = 1$,
  such that $s(C) = e(C) = c$, $h^0(\sI_C(c)) = a_2$, $h^1(\sI_C(c)) = r$,
  $h^1(\sO_C(c)) = b_1$ and $a_1 = 0, b_2 = 0$. Hence
\begin{equation*}
  {_0\!\hom_R}(I , M) = ra_2 \neq 0 \ \ {\rm and} \ \
  {_0\!\hom_R}(M , E) = rb_1 \neq 0
\end{equation*}
by \eqref{cor61} and \eqref{cor63}. Since $a_2=\beta_{1,c}$ and
$b_1=\beta_{2,c+4}$ the curves are obstructed by Theorem~\ref{mainthm}. To
see the existence, put $a = a_2$ and $b = b_1$. If $a = 1$, we consider curves
with $\Omega$-resolution
\begin{equation*} 
 0 \rightarrow \sO_{\proj{}}(-2)^{3r-1} \oplus
    \sO_{\proj{}}(-4)^{b} \rightarrow \sO_{\proj{}} 
    \oplus \Omega^{r}
    \oplus \sO_{\proj{}}(-3)^{b-1} \rightarrow \sI_C(c)
    \rightarrow 0
\end{equation*}
By Chang's results (\cite{C} or \cite{W1}, Thm. 4.1) there exists smooth
connected curves having $\Omega$-resolution as above. Moreover $c = 1 + b +
2r$, the degree $d = { c+4 \choose 2} - 3r - 7$ and the genus $g = (c+1)d - {
  c+4 \choose 3}+ 5$. If $a > 1$, curves with $\Omega$-resolution
\begin{equation*} 
  0 \rightarrow \sO_{\proj{}}(-1)^{a-2} \oplus
    \sO_{\proj{}}(-2)^{3r} \oplus 
    \sO_{\proj{}}(-4)^{b} \rightarrow  \sO_{\proj{}}^{a} \oplus
    \Omega^{r} 
    \oplus \sO_{\proj{}}(-3)^{b-1} \rightarrow \sI_C(c) \rightarrow 0
\end{equation*}
exist, they are smooth and connected (\cite{C} or \cite{W1}, Thm. 4.1), $c = a
+ b + 2r + 1$, $d = { c+4 \choose 2} - 3a - 3r - 6 $ and the genus $g = (c+1)d
- { c+4 \choose 3} + 3a + 3$. We leave the verification of details to the
reader, recalling only the exact sequences we frequently used in the
verification;
\begin{equation} \label{ex11}
 0 \rightarrow \Omega \rightarrow \sO_{\proj{}}(-1)^{4}
  \rightarrow \sO_{\proj{}} \rightarrow 0\ \ {\rm and} \ \ 0 \rightarrow
  \sO_{\proj{}}(-4) 
  \rightarrow \sO_{\proj{}}(-3)^{4} \rightarrow
  \sO_{\proj{}}(-2)^{6} \rightarrow  \Omega \rightarrow 0 
\end{equation}
Putting the two sequences together, we get the Koszul resolution of the
regular sequence $\{ X_0,X_1,X_2,X_3\} $.

We will analyze these curves a little further, using Laudal's description of
the completion of $O_{\HH_{(d,g)},(C)}$ (\cite{L1}, Thm.\! 4.2.4). This
completion is $k[[\HH^0(\sN_C)^{\vee}]]/o(\HH^1(N_C)^{\vee})$, where $o$ is a
certain obstruction morphism (giving essentially the cup and Massey products).
Now, consulting for instance the proof of Proposition~\ref{cupprod2}, we see
that the dual spaces of ${_0\!\Hom_R}(I,M)^{\vee}$ and
${_0\!\Hom_R}(M,E)^{\vee}$ inject into $\HH^0(\sN_C)^{\vee}$ and their
intersection is empty. This implies
$$
 \HH^0(N_C)^{\vee} \cong  T_{\gamma,\rho}^{\vee} \oplus {_0\Hom_R}(I, M)^{\vee}
 \oplus {_0\Hom_R}(M , E)^{\vee} \ \ {\rm as \ k-vector spaces}, $$
 and we can represent $k[[\HH^0(\sN_C)^{\vee}]]$ as
 $k[[Y_1,..Y_m,Z_{11},..,Z_{ar},W_{11},..,W_{rb}]]$, letting $Y_1,..Y_m$, resp.
 $Z_{11},..,Z_{ar}$, resp. $W_{11},..,W_{rb}$ correspond to a basis of
 $T_{\gamma,\rho}^{\vee}$, resp. ${_0\!\Hom_R}(I, M)^{\vee}$, 
 resp. ${_0\!\Hom_R}(M, E)^{\vee}$. Since $a_1 = 0, b_2 = 0$,
 we get by \eqref{cor61} and \eqref{cor63};
 \begin{equation*}
  {_{-4}\!\Hom_R}(I , M) = 0 \ \ {\rm and} \ \  {_{-4}\!\Hom_R}(M , E)
  = 0 \ .
\end{equation*}
By Remark~\ref{remunobstr} and Definition~\ref{delta}, $h^1(\sN_C) =
\delta^2(0) = a_2b_1$, and we can use Proposition~\ref{cupprod2} and its proof
to conclude that, modulo $\mathfrak m_\sO^3$ ($\mathfrak m_\sO$ the maximal
    ideal of the completion of $\sO_{\HH(d,g),(C)}$), we have
\begin{equation} \label{ex12}
 \sO_{H(d,g),(C)}/{\mathfrak m_{\sO}^3} =
       k[[Y_1,..Y_l,Z_{11},..,Z_{ar},W_{11},..,W_{rb}]]/{\mathfrak a} 
\end{equation}
where the ideal ${\mathfrak a}$ is generated by the components of the matrix
given by the product
\begin{equation}
  \label{ex13}
  \left[
    \begin{matrix}
      Z_{11} & \dots & Z_{1r} \\
      Z_{21} & \dots & Z_{2r} \\
      \vdots && \vdots \\
      Z_{a1} & \dots & Z_{ar}
    \end{matrix}
  \right]
  \left[
    \begin{matrix}
      W_{11} & \dots & W_{1b} \\
      W_{21} & \dots & W_{2b} \\
      \vdots && \vdots \\
      W_{r1} & \dots & W_{rb}
    \end{matrix}
  \right]
\end{equation}
Note that \eqref{ex13} corresponds precisely to the composition given by the
pairing of Proposition~\ref{cupprod1}! As in \cite{W1}, proof of Thm.\! 0.5, we
believe that the Massey products corresponding to \eqref{ex13} vanish, i.e.
the right-hand side of \eqref{ex12} is exactly the completion of
$\sO_{\HH(d,g),(C)}$.

The simplest case is $(r,a_2,b_1) = (1,1,1)$, which yields curves $C$ with
$s(C) = 4, d= 18$ and $g = 39$ (Sernesi's example \cite{Se} or \cite{EF}),
while the case $(r,a_2,b_1) = (2,1,1)$ yields curves $C$ with $s(C) = 6, d =
32$ and $g = 109$. More generally, the curves of the case $(r,1,1)$ satisfy
$h^1(\sN_C) = a_2b_1 = 1$, i.e. the ideal ${\mathfrak a}$ of \eqref{ex12} is
generated by the single element
\begin{equation} \label{ex14}
\sum_{i\, =\, 1 }^{r} ~ Z _{1i} \cdot W _{i\, 1} 
\end{equation}
For Sernesi's example $(r = 1)$, we recognize the known fact that this curve
sits in the intersection of two irreducible components of $\HH(d,g)$, while
for $r > 1$, the irreducibility of \eqref{ex14} can be used to see that $C$
belongs to a unique irreducible component of $\HH(d,g)$. Other examples of
singularities of $\HH(d,g)$ which belong to a unique irreducible component are
known (\cite{K1}, Rem. 3b\rm) \it{and \cite{Gu}, Thm. 3.10). In the next
  section we prove the irreducibility/reducibility by studying in detail the
  possible generizations of a Buchsbaum curve}.
\end{example} 

\section {The minimal resolution of a general space curve}
In this section we study generizations of space curves $C$ and how suitable
generizations will simplify the minimal resolution of $I(C)$. By a {\it
  generization} we mean a deformation to a ``more general curve'', cf.
Subsection 1.1. The general philosophy is that a sufficiently general curve of
any irreducible component of $\HH(d,g)$ should have as few repeated direct
factors "as possible" in consecutive terms of the minimal resolution. We prove
below a general result in this direction (Theorem~\ref{mainres}) and a more
restricted one (Proposition~\ref{mainres2}) for curves with special Rao
modules, using some nice ideas from \cite{MDP1} where they make explicit some
cancellations in the minimal resolution under flat deformation, in a special
case ($M \cong k$) which has the potential of being generalized. More
recently several papers have appeared using ``consecutive cancellations'' to
relate graded Betti numbers with the same Hilbert function (see \cite{P},
\cite{MI} and its references). Recalling the notations \eqref{resoluM}
and \eqref{resoluMI} from Rao's theorem (\cite{R}, Thm.\! 2.5), we show

\begin{theorem} \label{mainres}
  Let $C$ be a curve in $\proj{3}$ with postulation $\gamma$ and Rao module
    $M=M(C)$ and suppose the homogeneous ideal $I(C)$ has a minimal free
    resolution of graded $R$-modules;
\begin{equation}  \label{resoluMI2}
  0 \rightarrow L_4 \stackrel{\sigma \oplus 0}{\longrightarrow} L_3 \oplus F_2
    \rightarrow F_1 \rightarrow I(C) \rightarrow 0 \ . \ 
\end{equation}
If there exists a direct free factor $F$ satisfying $F_2 \cong F_2' \oplus F$
and $F_1 \cong F_1' \oplus F$, then there is a generization $C' \subseteq
\proj{3}$ of $C \subseteq \proj{3}$ in the Hilbert scheme $\HH(d,g)$ (in fact
in $\HH_{\gamma, M}$, i.e. with constant postulation and Rao module) whose
homogeneous ideal $I(C')$ has a minimal free resolution of the following form
\begin{equation*} 
  0 \rightarrow L_4 \stackrel{\sigma \oplus 0}{\longrightarrow} L_3 \oplus F_2'
    \rightarrow F_1' \rightarrow I(C') \rightarrow 0 \ . \ 
\end{equation*}
\end{theorem}

Now suppose $M=M(C)$ admits an $R$-{\it module} decomposition
$M=M' \oplus M_{[t]}$ where the diameter of $ M_{[t]}$ is $1$ (e.g. $C$ is
Buchsbaum). Let $ 0 \rightarrow L_4'  \stackrel{\sigma'}{\longrightarrow}
L_3' \rightarrow L_2' \rightarrow L_1' \rightarrow L_0' \rightarrow M'
\rightarrow 0$ be the minimal resolution of $M'$ and let $$
0 \rightarrow
R(-t-4)^r  \xrightarrow{\sigma_{[t]}} R(-t-3)^{4r} \rightarrow ... \rightarrow
R(-t)^r \rightarrow M_{[t]} \rightarrow 0$$
be the corresponding resolution of
$ M_{[t]}$ (which is ``$r$ times'' the Koszul resolution of the $R$-module $k
\cong R/(X_0, X_1, X_2, X_3)$.) By the Horseshoe lemma the minimal resolution
of $M$ is the direct sum of these two resolutions.  Looking to
\eqref{bettiMI}, we get $a_1 \cdot b_1 = 0$ and $a_2 \cdot b_2$ = 0 for a {\rm
  general} curve $C$ of $\HH(d,g)$ by Theorem~\ref{mainres}. Hence the
corresponding singularities of $\HH(d,g)$ given by Corollary~\ref{mcor} can
not occur for a general $C$, neither can the remaining class of singularities
due to

\begin{proposition} \label{mainres2} 
  Let $C$ be a curve in $\proj{3}$ and let $M(C) \cong M' \oplus M_{[t]}$ as
  $R$-{\it modules} where $ M_{[t]}$ is $r$-dimensional of diameter $1$ and
  supported in degree $t$. Moreover suppose the homogeneous ideal $I(C)$ has a
  minimal resolution of the following form;
\begin{equation} \label{resoluMI4}
  0 \rightarrow L_4' \oplus R(-t-4)^r  \xrightarrow{\sigma' \oplus
  \sigma_{[t]} \oplus 0} L_3' \oplus
  R(-t-3)^{4r} \oplus 
  F_2  \rightarrow F_1 \rightarrow I(C) \rightarrow 0 \ , \ 
\end{equation}
where $F_2 \cong P_2' \oplus R(-t-4)^{b _1}$ 
 and $F_1 \cong P_1' 
 \oplus R(-t)^{a_2}$ and where $P_2'$ (resp.  $P_1'$) is without direct free
 factors generated in degree $t+4$ (resp. $t$).
 
 $(a)$ \ \ Let $r \cdot b_1 \neq 0$ and let $m_1$ be a number satisfying $0
 \leq m_1 \leq \min \{r,b_1\}$. Then there is a generization $C' \subseteq
 \proj{3}$ of $C \subseteq \proj{3}$ in $\HH(d,g)$ (in fact in $\HH_{\gamma}$,
 i.e. with constant postulation $\gamma$) such that $I(C')$ has a
 free resolution of the following form;
\begin{equation*} 
 0 \rightarrow L_4' \oplus R(-t-4)^{r-m_1} 
    \rightarrow L_3' \oplus R(-t-3)^{4r} \oplus  P_2' \oplus R(-t-4)^{b_1-m_1}
    \rightarrow F_1 \rightarrow I(C') \rightarrow 0 \ , \ 
\end{equation*}
and such that $M(C') \cong M' \oplus M(C')_{[t]}$ as $R$-modules for some
$r-m_1$ dimensional module $M(C')_{[t]}$ supported in degree $t$. The
resolution is minimal except possibly in degree $t+3$ where some of the common
free factors of $R(-t-3)^{4r}$ and $F_1$ may cancel. Moreover if $L_2'$ does
not contain a direct free factor generated in degree $t+4$, then $ \
{_0\!\hom_R}(M(C')_{[t]},E(C'))=(r-m_1)(b_1-m_1)$.

$(b)$ \ Suppose $L_2'$ is without direct free factors generated in
degree $t$. If $r \cdot a_2 \neq 0$ and if $m_2$ is a number satisfying $0
\leq m_2 \leq \min \{r,a_2\}$, then there is a generization $C' \subseteq
\proj{3}$ of $C \subseteq \proj{3}$ in $\HH(d,g)$ (with constant
specialization) such that $I(C')$ has a minimal free resolution of the
following form;
\begin{equation*} 
 0 \rightarrow L_4' \oplus R(-t-4)^{r-m_2} 
    \rightarrow L_3' \oplus G_2 \rightarrow G_1 \oplus
    R(-t)^{a_2-m_2}  \rightarrow I(C') \rightarrow 0 \  
\end{equation*}
for some $R$-free modules $G_2$ and $G_1$ where $G_1$ is without direct free
factors generated in degree $t$. Moreover $M(C') \cong M' \oplus M(C')_{[t]}$
as $R$-modules for some $r-m_2$ dimensional module $M(C')_{[t]}$ supported in
degree $t$, and we have $ \ {_0\!\hom_R}(I(C'),
M(C')_{[t]})=(r-m_2)(a_2-m_2)$.
\end{proposition}

Once we have proved a key lemma, the proof of Theorem~\ref{mainres} is
straightforward while the proof of Proposition~\ref{mainres2} is a little bit
more technical. Note that the assumptions on $L_2'$ in
Proposition~\ref{mainres2}$(a)$ and $(b)$ show that ${_0\!\Ext_R^2}(M',M')=0
\Rightarrow {_0\!\Ext_R^2}(M,M)=0$ (Remark~\ref{remmcor}), indicating that our
results of this section combine nicely with Theorem~\ref{obstr}. We delay the
proof of these results until the end of this section.

Now combining these two results with Theorem~\ref{mainthm} in the diameter one
case, we get

\begin{corollary} \label{mainrescor} 
  Let $C$ be a curve in $ \proj{3}$ whose Rao module $M \neq 0$ is of diameter
  $1$ and concentrated in degree $c$, and let $\beta_{1,c+4}$ and
  $\beta_{1,c}$ (resp.  $\beta_{2,c+4}$ and $\beta_{2,c}$) be the number of
  minimal generators (resp. minimal relations) of degree $c+4$ and
  $c$ respectively. 
  
  (a) If $C$ is generic in $\HH_{\gamma, \rho}$, then $\HH_{\gamma}$ is 
  smooth at $(C)$. Moreover $C$ is obstructed if and
  only if $\beta_{1,c} \cdot \beta_{2,c+4} \neq 0$. Furthermore if
  $\beta_{1,c} =0$ and $ \beta_{2,c+4}= 0$, then  $C$ is generic in
  $\HH(d,g)$. 
  
  (b) If $C$ is generic in $\HH_{\gamma}$, then $C$ is unobstructed. Indeed
  both $\HH(d,g)$ and $\HH_{\gamma}$ are smooth at $(C)$.  In particular every
  irreducible component of $\HH(d,g)$ whose generic curve $C$ satisfies $\diam
  M(C) \leq 1$ is {\it reduced} (i.e. generically smooth).
\end{corollary}

\begin{proof} $(a)$ $C$ is generic in $\HH(d,g)$ by  
  Proposition~\ref{gen} because ${_0\!\Hom_R}(M,E)={_0\!\Hom_R}(I,M)=0$ by
  \eqref{cor61} and \eqref{cor63}. The other statements follow directly from
  Theorem~\ref{mainthm}, Theorem~\ref{mainres} and
  Proposition~\ref{remmainthm}.
  
  $(b)$ If $C$ is generic in $\HH_{\gamma}$, then we immediately have
  $\beta_{1,c} \cdot \beta_{2,c}= 0$ and $r \cdot \beta_{2,c+4} = 0$ by
  Theorem~\ref{mainres} and Proposition~\ref{mainres2}. Since $r > 0$ we see
  by Theorem~\ref{mainthm} that $\HH(d,g)$ (and of course $\HH_{\gamma}$ by
  $(a)$) is smooth at $(C)$.  Finally if $C$ is a generic curve of some
  irreducible component of $\HH(d,g)$ satisfying $\diam M(C) \leq 1$ and
  $\gamma$ is the postulation of $C$, then $C$ is generic in $\HH_{\gamma}$
  and we conclude easily.
\end{proof}

Corollary~\ref{mainrescor}$(a)$ generalizes \cite{BKM} Prop. 1.1 which tells
that a curve $C$ of maximal rank or maximal corank of $\diam M(C) = 1$, which
is generic in $\HH_{\gamma,\rho}$, is unobstructed.
  
Even though we can extend the next corollary to Buchsbaum curves satisfying
${_0\!\Ext_R^2}(M, M) = 0$ (i.e. ${_0\!\Hom_R}(L_2, M) = 0$), we have chosen
to formulate it for the somewhat more natural set of Buchsbaum curves $C$ of
$\diam M(C) \leq 2$. Note that Buchsbaum curves of maximal rank satisfy $\diam
M(C) \leq 2$ (\cite{MIG}, Cor.\! 3.1.4, \cite{EF1}, Cor.\! 2.8), and
Corollary~\ref{mainrescor} and \ref{mainrescor2} (and \cite{Mir}) give answers
to the problems on unobstructedness of Buchsbaum curves raised by Ellia and
Fiorentini in \cite{EF1}.

\begin{corollary}  \label{mainrescor2} 
  Let $C$ be a Buchsbaum curve of $\diam M(C) \leq 2$.
  Then there exists a generization $C'$ of $C$ in $H(d,g)$ such that $C'$ is
  Buchsbaum (or ACM with $L_4=0$) and such that the modules of the three sets
\begin{equation*}
\{F_2, F_1 \} \ ,\ \{ L_4, F_2 \} \ \ {\it and\/} \ \ \{  L_4, F_1(-4) \}
\end{equation*}
in its minimal resolution, $ 0 \rightarrow L_4 \stackrel{\sigma \oplus
  0}{\longrightarrow} L_3 \oplus F_2 \rightarrow F_1 \rightarrow I(C')
\rightarrow 0,$ are without common direct free factors. Hence
${_0\!\Hom_R}(I(C'),M(C')) = {_0\!\Hom_R}(M(C'), E(C'))= 0$ and $\HH(d,g)$ is
smooth at $(C')$.
\end{corollary}
  
\begin{proof}  Firstly note that since the module structure of $M$ of any
  Buchsbaum curve is trivial, we get from the resolution \eqref{resoluMI2}
  that $ {_0\!\Hom_R}(I , M) \cong {_0\!\Hom_R}(F_1, M)$.  Since $M \cong \ker
  \HH_{*}^3(\tilde {\sigma} \oplus 0)$, it follows that the latter group
  vanishes if and only if $L_4$ and $F_1(-4)$ are without common direct free
  factors. Moreover by arguing as in the proof of Corollary~\ref{mcor} we get
  $ {_{-4}\!\Hom_R}(F_2 , M)^{\vee} \cong {_0\!\Hom_R}(M , E)$ which vanishes
  if and only if $L_4$ and $F_2$ are without common direct free factors.
  
  Now, by Theorem~\ref{mainres}, $ \{F_2, F_1 \}$ have no common direct free
  factors, and writing $M(C) \cong M_{[c]} \oplus M_{[c-1]}$ as $R$-modules,
  we can successively apply Proposition~\ref{mainres2} to $M_{[c]}$ and
  $M_{[c-1]}$. Indeed the $(a)$ part of Proposition~\ref{mainres2} with
  $M_{[t]}=M_{[c]} $ and $ m_1 = \min \{r,b_1\}$ shows that $\{ L_4, F_2 \}$
  for some generization of $C$ are without common direct free factors of
  degree $c+4$.  Then we proceed by $(b)$ to see that $\{ L_4, F_1(-4) \}$ for
  some further generization of $C$ are without common direct free factors of
  degree $c+4$.  Similarly we use Proposition~\ref{mainres2} with
  $M_{[t]}=M_{[c-1]} $ to see that there remains, up to a suitable
  generization $C'$, also no common direct free factor of degree $c+3$ in $\{
  L_4, F_2 \} $ and $ \{ L_4, F_1(-4) \}$.  Hence we have
  ${_0\!\Hom_R}(I(C'),M(C')) = {_0\!\Hom_R}(M(C'), E(C'))= 0$ by the first
  part of the proof and we conclude by Proposition~\ref{gen}.
\end{proof}

We should have liked to generalize Corollary~\ref{mainrescor2} to the arbitrary
case of diameter 2 by dropping the Buchsbaum assumption. In particular if we
could prove a result analogous to Corollary~\ref{mainrescor2} for curves whose
Rao module $M$ is the generic module of diameter two (cf.  \cite{MDP2} for
existence and minimal resolution), we would be able to answer affirmatively
the following question (which we believe is true). \\

{\it Question.\/} Is any irreducible component of $\HH(d,g)$ whose Rao module
of its generic curve is concentrated in at most two consecutive degrees,
generically smooth? \\

In our corollaries we have used Theorem~\ref{mainres} and
Proposition~\ref{mainres2} to consider generic curves, or to get the existence
of a certain generization, with nice obstruction properties.  We may, however,
also use our results to study many different generizations of a given curve
$C$, see the works of Amasaki, Ellia and Fiorentini and others (\cite{A},
\cite{Se}, \cite{EF}, \cite{K3}) for similar approaches. Hence we may see when
$C$ sits in the intersection of different integral components of $\HH(d,g)$.
There may be quite a lot of such irreducible components of $\HH(d,g)$
\cite{Gi}. We will soon look closely to the possible generizations of a curve
of diameter one in the case $\beta_{1,c} \cdot \beta_{2,c+4} \neq 0$.  To get
a flavour of the other possibilities, we consider the following example of a
non-generic curve of $ \HH_{\gamma,M}$.

\begin{example} \label{ex311}
  In \cite{BKM} and \cite{W1} one proves the existence of an obstructed
  curve of $\HH(33,117)_S$ of maximal rank with one-dimensional Rao module.
  Since the degrees of the minimal generators of $I(C)$ are given in
  \cite{BKM} and $M = H^1(\sI_C(5))$, we easily find the minimal resolution to
  be
\begin{equation*}
0 \rightarrow R(-9) \rightarrow R(-10)^{2} \oplus  R(-9)
\oplus  R(-8)^{4} \rightarrow R(-9) \oplus  R(-8) \oplus
 R(-7)^{5} \rightarrow I(C) \rightarrow 0 \ .
\end{equation*}
It follows from Theorem~\ref{mainthm} of this paper that $C$ is obstructed. By
Proposition~\ref{mainres2} (resp. Theorem~\ref{mainres}) there exists a
generization $C_1$ (resp. $C_2$) of $C$, obtained by removing the direct
factor $R(-9)$ from $L_4$ and $F_2$ (resp. from $F_2$ and $F_1$). The curve
$C_1$ is ACM, hence unobstructed, and belongs to a unique irreducible
component $V$ of $\HH(33,117)_S$. Moreover the curve $C_2$ is unobstructed by
Theorem~\ref{mainthm}. Now looking only to the semicontinuity of
$h^1(\sI_C(5))$ and $h^1(\sO_C(5))$, there is a priori a possibility that
$C_2$ may belong to $V$. By Corollary~\ref{mainrescor}$(a)$ or by
Proposition~\ref{gen}, however, $C_2$ is generic in $\HH(33,117)_S$ since we
may suppose $C_2$ is generic in $\HH(33,117)_{\gamma}$. Hence the irreducible
component $W$ of $\HH(33,117)_S$ to which $C_2$ belongs, satisfies $W \neq
V$!!  Since $C$ is contained in the intersection of the components, we get the
main example of \cite{BKM} from our results.
\end{example}

As an illustration of the main results of this section, we restrict to curves
which are generic in $ \HH_{\gamma,M}$, or more generally to curves which
satisfy $a_1 \cdot b_1 = 0$ and $a_2 \cdot b_2 = 0$ (letting $a_1=
\beta_{1,c+4}$, $a_2=\beta_{1,c}$, $b_1= \beta_{2,c+4}$ and $b_2=
\beta_{2,c}$). Thus we consider the case
\begin{equation}  \label{exfi1}
   a_1 = 0 ,\  b_2 = 0 \ \ {\rm and} \ \ (  a_2 \neq 0 \ {\rm or} \  b_1 \neq
   0 ) 
\end{equation}
where proper generizations as in Proposition~\ref{mainres2} occur, to give a
rather complete picture of the existing generizations in $\HH(d,g)$ (caused by
simplifications of the minimal resolution). Let $n(C) = (r,a_1,a_2,b_1,b_2)$ be
an associated $5$-tuple. Only for curves satisfying $a_1 = 0$ and $b_2 = 0$ we
allow the writing $n(C) = (r,a_2,b_1)$ as a triple. Thanks to \cite{B} we
remark that any curve $D$ satisfying $n(D) = n(C)$ and $\gamma_{D}(v) =
\gamma_{C}(v)$ for $v \neq c$, belongs to the same {\it irreducible} family $
\HH_{\gamma,M}$ as $C$, i.e. a further generization of $C$ and $D$ in $
\HH_{\gamma,M}$ lead to the "same" generic curve. Now given a curve $C$ with
$n(C) = (r,a_2,b_1)$, we have by Proposition~\ref{mainres2}:
\begin{equation} \label{exfi2}
  \begin{gathered}
    \text{For any pair $(i,j)$ of non-negative integers such that }
    \text{$r-i-j \geq 0$, $a_2-i \geq 0$ and $b_1-j \geq 0$,} \\
    \text{there exists a generization $C_{ij}$ of $C$ in $\HH(d,g)$ such that
      $n(C_{ij}) = (r-i-j,a_2-i,b_1-j)$}.
  \end{gathered}
\end{equation}  

Note that if we link $C$ to $C_l$ as in Proposition~\ref{mlink}, we get, by
combining \eqref{thm66}, \eqref{cor61} and \eqref{cor63} that the $5$-tuple
$n(C_l) = (r(C_l),a_1(C_l),a_2(C_l),b_1(C_l),b_2(C_l))$ is equal to
$(r,b_2,b_1,a_2,a_1)$ where $n(C) = (r,a_1,a_2,b_1,b_2)$. In particular if $C$
satisfies \eqref{exfi1}, then the linked curve $C_l$ also does.

As an example, let $n(C) = (4,3,2)$ (such curves exist by Example~\ref{ex1}).
By \eqref{exfi2} we have $10$ different generizations $C_{ij}$ among which two
curves correspond to the triples $n(C_{22}) = (0,1,0)$ and $n(C_{31}) =
(0,0,1)$, i.e. they correspond to two {\it unobstructed} ACM curves with
different postulation. Hence they belong to two different irreducible
components of $\HH(d,g)$ having (C) in their intersection. Pushing this
argument further, we get at least

\begin{proposition}  \label{propfi}
  Let $C$ be a curve in $\proj{3}$ whose Rao module $M \neq 0$ is
  $r$-dimensional and concentrated in degree $c$, let $a_1= \beta_{1,c+4}$ and
  $a_2=\beta_{1,c}$ (resp.  $b_1= \beta_{2,c+4}$ and $b_2= \beta_{2,c}$) be
  the number of minimal generators (resp. minimal relations) of degree $c+4$
  and $c$ respectively, and suppose
\begin{equation*}
  a_1 = 0 , \ b_2 = 0 \ \ \ {\rm and}  \ \ \  a_2 \cdot b_1 \neq 0 \ .
\end{equation*}

(a) If \ $r < a_2 + b_1$, then $C$ sits in the intersection of at least two
irreducible components of $\HH(d,g)$. Moreover, the generic curve of any
component containing $C$ is arithmetically Cohen-Macaulay, and the number
$n(comp,C)$ of irreducible components containing $C$ satisfies
\begin{equation*}
  \min\{a_2,r\} + min\{b_1,r\} - r + 1 \leq n(comp,C) \leq r + 1 \ .
\end{equation*}
 In the case $s(C) = e(C) = c$, we have equality to the left.
 
(b) If \ $r \geq a_2 + b_1$ and $s(C) = e(C) = c$, then $C$ is an obstructed
 curve which belongs to a unique irreducible component of $\HH(d,g)$.
\end{proposition}

\begin{proof} We firstly prove $(b)$. Let $C'$ be any generization of $C$ in
  $\HH(d,g)$ and let $n(C') = (r',a_1',a_2',b_1',b_2')$ be the associated
  5-tuple where $r'=0$ corresponds to the ACM case of $C'$. Since $s(C) = c$
  and since the number $s(C)$ increases under generization by the
  semicontinuity of $h^0(\sI_C(v))$, we get $s(C') \geq c$ as well as
  $h^0(\sI_{C'}(c)) = a_2'$ and $b_2' = 0$.  Similarly $e(C) = c$
  implies $h^1(\sO_{C'}(c)) = b_1'$ and $a_1' = 0$.  Applying these
  considerations to $C' = C$, we get $\chi(\sI_{C}(c)) \leq 0$ by the
  assumption $r \geq a_2 + b_1$.
  
  Now let $C'$ be the generic curve of an irreducible component containing
  $C$. By Proposition~\ref{mainres2} we get $r'a_2' = 0$ and $r'b_1' = 0$
  which combined with $\chi(\sI_{C'}(c)) = \chi(\sI_{C}(c)) \leq 0$ yields
  $a_2' = 0$ and $b_1' = 0$.  Hence $n(C') = (r-a_2-b_1,0,0,0,0)$ for any
  generic curve of $\HH(d,g)$. Since $\gamma_{C'}(v) = \gamma_{C}(v)$ for $v
  \neq c$ by semicontinuity and the vanishing of $\HH^1(\sI_C(v))$, any such
  $C'$ belongs to the same irreducible component of $\HH(d,g)$ by the
  irreducibility of $\HH_{\gamma_{C'},M(C')}$. Moreover
  $C$ is obstructed by Theorem~\ref{mainthm}, and $(b)$ is proved.
  
  $(a)$ Suppose $ r < a_2 + b_1$. To get the lower bound of $n(comp,C)$ (which
  in fact is $\geq 2$), we use \eqref{exfi2} to produce several generic curves
  of $\HH(d,g)$ which are generizations of $C$. Indeed let $m(a) =
  \min\{a_2,r\}$ and $m(b) = min\{b_1,r\}$. By \eqref{exfi2} there exist
  generizations $C_0$, $C_1$,..,$C_{m(a)+m(b)-r}$ such that $n(C_0) =
  (0,a_2-m(a),b_1+m(a)-r)$, $n(C_1) = (0,a_2-m(a)+1,b_1+m(a)-r-1)$,...,
  $n(C_{m(a)+m(b)-r}) = (0,a_2+m(b)-r,b_1-m(b))$. Since the curves $C_i$ are
  ACM and have different postulations, they belong to $m(a) + m(b) - r + 1$
  different components, and we get the minimum number of irreducible
  components as stated in the proposition.
  
  To see that the generic curve $C'$ of {\it any\/} component containing $C$
  is ACM, we recall that $r'a_2' = 0$ and $r'b_1' = 0$ by
  Proposition~\ref{mainres2} with notations as in the first part of the proof.
  Suppose $r' \neq 0$. Then $a_2' = 0$ and $b_1' = 0$. To get a contradiction,
  we remark that $\gamma_{C'}(v) = \gamma_{C}(v)$ for $v <
  c$, 
  from which we get $h^0(\sI_{C'}(c)) + b_2' = h^0(\sI_{C}(c)) - a_2 $ since
  $a_2$ (resp. $b_2'$) is the only possibly non-vanishing graded Betti number
  of $I(C)$ (resp. $I(C')$) in degree $c$. Hence $h^0(\sI_{C'}(c)) \leq
  h^0(\sI_{C}(c)) - a_2$ and similarly we have the "dual" result
  $h^1(\sO_{C'}(c)) \leq h^1(\sO_{C}(c)) - b_1$. Adding the inequalities, we
  get
\begin{equation*}
\chi(\sI_{C'}(c)) + h^1(\sI_{C'}(c))
 \leq \chi(\sI_{C}(c)) + h^1(\sI_{C}(c)) - a_2 - b_1 < \chi(\sI_{C}(c))\ ,  
\end{equation*}
i.e. a contradiction because $\chi(\sI_{C'}(c)) = \chi(\sI_{C}(c))$.  Now
using the fact that the generic curve $C'$ of any irreducible component
containing $C$ is ACM and that $\HH_{\gamma_{C'},M(C')}$
is irreducible, we prove easily that $n(comp,C) \leq r + 1$ because there are
at most $r + 1$ different postulations $\gamma_{C'}$. Indeed since $M(C') = 0$,
$\gamma_{C'}(v) = \gamma_{C}(v)$ for $v \neq c$ and
\begin{equation*}
  \gamma_{C'}(c) + \sigma_{C'}(c) = \chi(\sI_{C'}(c)) = \chi(\sI_{C}(c))=
      \gamma_{C}(c) + \sigma_{C}(c) - r
\end{equation*}
where $\sigma_C(v) = h^1(\sO_C(v))$, we see that the different choices of
$\gamma_{C'}$ can happen in degree $v = c$ only, and that they are given by
$\gamma_{C'}(c) = \gamma_C(c)-i$ where $i$ is chosen among $\{0,1,2,..,r\}$.

Suppose $s(C) = e(C) = c$. Since in this case $\gamma_C(c) = a_2$
and $\sigma_C(c) = b_1$ by arguments as in the first part of the proof, we
can easily limit the (at most) $r+1$ different choices of the postulation
$\gamma_{C'}(c) = \gamma_C(c)-i$ above by choosing
$$
   m(a) \leq i \leq r - m(b)  $$
i.e. $n(comp,C)$ equals precisely  $m(a) + m(b) - r + 1$, and we are done.
\end{proof}

\begin{example}  Now we reconsider some particular cases of Example~\ref{ex1},
  even though Proposition~\ref{propfi} is well adapted to treat the whole
  example in detail. Recall that for any triple $(r,a_2,b_1)$ of natural
  numbers, there exists a smooth connected curve $C$ with $n(C) = (r,a_2,b_1)$
  and $s(C) = e(C) = c(C)$ by Example~\ref{ex1}. In particular
  
  $(a)$ For every integer $r > 0$ there exists a smooth connected curve $C$,
  with triple $n(C) = (r,r,r)$, of degree $d$ and genus $g$ as in
  Example~\ref{ex1}, which is contained in $r + 1$ irreducible components of
  $\HH(d,g)_S$.  Moreover the generic curves of all the components containing
  $C$ are ACM.
  
  $(b)$ For every $r > 0$ there exists an obstructed, smooth connected curve
  with triple $(r,a_2,b_1) = (2t,t,t)$ or $(2t+1,t,t)$, of degree $d$ and
  genus $g$ as given by Example~\ref{ex1}, which belongs to a unique
  irreducible component of $\HH(d,g)_S$ by Proposition~\ref{propfi}. In
  particular the obstructed curve $C$ with $(r,a_2,b_1) = (2,1,1)$ belongs to
  a {\it unique\/} irreducible component of $\HH(32,109)_S$, confirming what
  we saw in Example~\ref{ex1}.
\end{example}

To prove Theorem~\ref{mainres} and Proposition~\ref{mainres2} we need a lemma
for deforming a module $N$, which basically is known (and related to
\cite{MDP1}, Prop.\! 2.1, p.\! 140). For our purpose it suffices to see that
if we can lift a (three term) resolution with augmentation $N$ to a {\it
  complex}, then the complex defines a flat deformation of $N$. In the case $N
= I(C)$ where $C$ has e.g. codimension $2$ in $\proj{3}$, we also know that a
deformation of an ideal $I(C)$ is again an ideal, i.e.

\begin{lemma} \label{lemgen} 
 Let $C$ be a curve in $\proj{3}$ whose homogeneous ideal $I(C)$
  has a minimal resolution of the following form 
\begin{equation*}
(L^{\bullet})\ \ \ \ \ \ \ \ \ \ \  \ 0 \rightarrow \bigoplus_i
R(-i)^{\beta_{3,i}} \stackrel{\varphi}{\longrightarrow} \bigoplus_i 
R(-i)^{\beta_{2,i}}  \stackrel{\psi}{\longrightarrow} \bigoplus_i
R(-i)^{\beta_{1,i}} 
\rightarrow I(C) \rightarrow 0 \ . \ \ \ \ \  
\end{equation*}
Let $A$ be a finitely generated $k$-algebra, $B$ the localization of $A$ in a
maximal ideal $\wp$, and suppose there exists a complex
\begin{equation*}
(L_B^{\bullet}) \ \ \ \  \ \ \ \  \bigoplus_i R_B(-i)^{\beta_{3,i}}
 \stackrel{\varphi_B}{\longrightarrow}  \bigoplus_i 
R_B(-i)^{\beta_{2,i}} \stackrel{\psi_B}{\longrightarrow} \bigoplus_i
 R_B(-i)^{\beta_{1,i}}  ~ ~,  ~ ~R_{B}  =  R \otimes_k B \ ,  
\end{equation*}
such that $L_B^{\bullet} \otimes_B (B/{ \wp}) \cong L^{\bullet}$. Then
$(L_B^{\bullet})$ is acyclic, $\varphi_B$ is injective and the cokernel of
$\psi_B$ is a flat deformation of $I(C)$ as an ideal (so $\coker(\psi_B)
\subseteq R_B$ defines a flat deformation of $C \subseteq \proj{3}$ with
constant postulation). Moreover for some $a \in A- \wp$, we can extend this
conclusion to $A_a$ via $\Spec(B) \hookrightarrow Spec(A_a)$, i.e. there
exists a flat family of curves $C_{\Spec(A_a)} \subseteq \proj{3} \times
\Spec(A_a)$ whose homogeneous ideal $I(C_{A_a})$ has a resolution (not
necessarily minimal) of the form
\begin{equation*}
(L_{A _{a}}^{\bullet}) ~ ~ ~  ~ ~ 0 \rightarrow \bigoplus \, R _{ A _{a}
  }(-i)^{\beta_{3,i}}\rightarrow  \bigoplus \, R _{A_{a}}(-i)^{\beta_{2,i}}
  \rightarrow \bigoplus \, R_{A_{a}}(-i)^{\beta_{1,i}} \rightarrow I(C_{A_{a}})
  \rightarrow  0  ~ . ~ 
\end{equation*}
\end{lemma}

\begin {proof}[Proof (sketch)] \ If $ E = \coker \varphi$ and $E_B = \coker
  \varphi_B$, then one proves easily that $E_B \otimes_B(B/\wp) = E$,
  $\Tor_1(E_B,B/\wp) = 0$ and that $\varphi_B$ is injective. By the local
  criterion of flatness, $E_B$ is a flat deformation of $E$. Letting $Q_B =
  \coker (E_B \rightarrow \oplus_i R_B(-i)^{\beta_{1,i}})$, we can argue as we
  did for $ E_B$ to see that $Q_B$ is a flat deformation of $I(C)$ and that
  $L_B^{\bullet}$ augmented by $Q_B$ is exact.
  
  To prove that $Q_B$ is an ideal in $R_B$, we can use the isomorphisms
  $\HH^{i-1}(\sN_C) \cong \Ext_{\sO_{\proj{}}}^i(\tilde I,\tilde I)$ for $i =
  1,2$, 
  interpreted via deformation theory and repeatedly applied to $B_{i+1}
  \rightarrow B_{i}$ for $i\geq 1$ $(B_i=B/\wp^i)$, to see that a
  deformation of the $\sO_{\proj{}}$-Module $\tilde I$ (such as $\tilde
  {Q_B}$) corresponds to a deformation of the curve $C$ in the usual way, i.e.
  via the cokernel of $\tilde i$: $\tilde {Q_B} \rightarrow \tilde {R_B}$. We
  get in particular a morphism $ \HH_{*}^0(\tilde i)$: $Q_B \rightarrow R_B$
  which proves what we want (one may give a direct proof using Hilbert-Burch
  theorem (cf.  \cite{MDP1}, page 37-38)).

  Finally we easily extend the morphism $i$ and any morphism of the resolution
  $L_B^{\bullet}$ to be defined over $A_{a'}$, for some $a' \in A- \wp$ (such
  that $L_{A_{a'}}^{\bullet}$ is a complex). By shrinking $\Spec A_{a'}$ to
  $\Spec A_a$, $a \in A-\wp$, we get the exactness of the complex and the
  flatness of $I(C_{A_a})$ because these properties are open.
\end{proof}

\begin{proof}[Proof (of Theorem~\ref{mainres})]\ 
  Suppose that $F$ has rank $s$ and consider the $s$ by $s$ submatrix
  $M(\psi)$ of $\psi$ in
\begin{equation*}
 0 \rightarrow L_4  \xrightarrow{\sigma \oplus 0 \oplus 0} L_3 \oplus F_2'
 \oplus \, F  \stackrel{\psi}{\longrightarrow} F_1' \oplus \, F \rightarrow
 I(C) \rightarrow 0  \   
\end{equation*}
which corresponds to $F \rightarrow F$. As in the "Lemma de
g\'{e}n\'{e}risation simplifiantes" (\cite{MDP1}, page 189), we can change the
$0's$ on the diagonal of $M(\psi)$ to some $\lambda_1,...,\lambda_s$ where the
$\lambda_i's$ are indeterminates of degree zero. Keeping $\sigma \oplus 0
\oplus 0$ unchanged, we still have a complex which by Lemma~\ref{lemgen}
implies the existence a flat family of curves over $ \Spec(A_a)$, $A =
k[\lambda_1,...,\lambda_s]$, for some $a \in A-(\lambda_1,...,\lambda_s)$. Let
$\lambda:=\Pi_{i=1}^s \lambda_i$ be the product. Since any curve $C'$ of the
family given by $\Spec(A_{\lambda a})$ has a resolution where $F$ is redundant
($F$, and only $F$, is missing in its {\it minimal} resolution), and since we
may still interpret the Rao module $M(C')$ as $\ker \HH_{*}^3(\tilde {\sigma}
\oplus 0 \oplus 0)$ with ${\sigma \oplus 0 \oplus 0}$ as above (so the whole
family given by $\Spec(A_a)$ has constant Rao modules), we conclude easily.
\end{proof}

  \begin{remark} Slightly extending the proof and using Bolondi's result on the
    irreducibility of $\HH_{\gamma,M}$ (\cite{B}), one may prove that set $U$
    of points $(C)$ of the scheme $\HH_{\gamma,M}$ whose modules $F_2$ and
    $F_1$ of the minimal resolution \eqref{resoluMI2} of $I(C)$ are without
    common direct free factors, form an {\rm open} (and non-empty if a curve
    with minimal resolution \eqref{resoluMI2} exists) irreducible subset of
    $H_{\gamma,M}$.
\end{remark}

\begin{proof}[Proof (of Proposition~\ref{mainres2})] $(a)$
  Since we have the assumption that $M \cong M' \oplus M_{[t]}$ as $R$-{\it
    modules}, the minimal resolution \eqref{resoluM} of $M$ is given as the
  direct sum of the resolution of $M'$ and the one of $ M_{[t]}$ which is
  ``$r$-times'' the Koszul resolution associated with the regular sequence
  $\{X_0, X_1, X_2, X_3 \}$. The matrix associated to $\sigma_{[t]}$ (resp.
  $\sigma = \sigma' \oplus \sigma_{[t]}$) will have the form
 \begin{equation}  \label{ex23}
  \left[
    \begin{matrix}
      \underline{X} & 0 & \dots & 0 \\
      0 & \underline{X} & \dots & 0 \\
      \vdots & \vdots & \ddots & \vdots \\
      0 & 0 & \dots & \underline{X}
    \end{matrix}
  \right]
 \ \ \ \ \ \ \left[
  \begin{matrix}
      \sigma' &  0 \\
      0  & \sigma_{[t]} \\
   \end{matrix}
  \right]
\end{equation}
where $\underline{X}$ is $(X_0,X_1,X_2,X_3)^T$ and each "row" in the left
matrix is a $4 \times r$ matrix, etc. Let $\eta_j$: $R(-t-4) \rightarrow L_4'
\oplus R(-t-4)^r$ be the map into the $j$-th direct factor of $ R(-t-4)^r$,
$1\leq j \leq r$, and let $\pi_i$ : $L_3' \oplus R(-t-3)^{4r} \oplus P_2'
\oplus R(-t-4)^{b _1} \rightarrow R(-t-4)$ be the projection onto the $i$-th
factor of $R(-t-4)^{b _1}$, $1\leq i \leq b_1$. Similar to what was observed
by Martin Deschamps and Perrin in the case $M \cong k$ (\cite{MDP1}, page 189)
we can change the $0$ component in the matrix of $\sigma \oplus 0$ which
corresponds to $\pi_i \eta_j$: $ R(-t-4) \rightarrow R(-t-4)$, to some
indeterminate of degree zero. To get a complex we need to change four columns
of the matrix $A$ associated to $L_3 \oplus F_2 \rightarrow F_1$ as follows.
Let $r_1:= \rank F_1$ and look to the column $(a_k)$, $1 \leq k \leq r_1$, of
$A$ which corresponds to the map $R(-t-4) \rightarrow F_1$ from the $i$-th
factor of $ R(-t-4)^{b _1}$. Put $a_k= \sum_{l=0}^3 \gamma_{k,l}^i X_l$ for
every $1 \leq k \leq r_1$. Since the resolution is minimal, such
$\gamma_{k,l}^i$ exist, but they are not necessarily unique. Since the column
of the matrix of $\sigma \oplus 0$ which corresponds to $\eta_j$ consists of
only $0$'s and $\underline X$ (cf. \eqref{ex23}) there are precisely four
columns $[H_{k,0}^j,H_{k,1}^j,H_{k,2}^j,H_{k,3}^j]$, $1 \leq k \leq r_1$, of
$A$ satisfying $\sum_{l=0}^3 H_{k,l}^j X_l=0$ for every $k$ which may
contribute to the composition $(\sigma \oplus 0) \eta_j$. Now if we change the
trivial map $\pi_1 \eta_1$ to the multiplication by an indeterminate
$\lambda_1$ and simultaneously change the four columns
$[H_{k,0}^1,H_{k,1}^1,H_{k,2}^1,H_{k,3}^1]$ of $A$ to $[H_{k,0}^1-
\gamma_{k,0}^1\lambda_1,H_{k,1}^1-
\gamma_{k,1}^1\lambda_1,H_{k,2}^1-\gamma_{k,2}^1\lambda_1,H_{k,3}^1-
\gamma_{k,3}^1\lambda_1]$, leaving the rest of $A$ unchanged, we still get
that \eqref{resoluMI4} defines a complex. We can proceed by simultaneously
changing the $0$ component of $\pi_2 \eta_2$ to $\lambda_2$ and the
corresponding four columns of the matrix $A$ as described above, etc. Put
$\lambda:=\Pi_{i=1}^{m_1} \lambda_i$. By Lemma~\ref{lemgen} we get a flat
irreducible family of curves $C'$ over
$\Spec(k[\lambda_1,...,\lambda_{m_1}]_{a})$, for some $a \in
A-(\lambda_1,...,\lambda_{m_1})$, having the same (not necessarily minimal)
resolution, hence the same postulation, as $C$. Since $\lambda$ is invertible
in $\Spec(k[\lambda_1,...,\lambda_{m_1}]_{\lambda \cdot a})$, we can remove
redundant factors of the resolution of $I(C')$ in this open set. Since $M(C')
\cong \ker \HH_{*}^3(\tilde {\sigma} \oplus 0 \oplus 0)$, we have a
generization $C'$ with properties as claimed in Proposition~\ref{mainres2}.
Note that since we have changed $4m_1$ columns of $A$ we may have changed some
zero entries of $A$ to non-zero constants, making the resolution non-minimal
in degree $t+3$ (only). Finally using Remark~\ref{remmcor} for $v=0$, the
assumption on $L_2'$ shows that \eqref{cor62} holds and hence we conclude by
the left formula of \eqref{cor63}.

$(b)$ \ We will prove $(b)$ by linking $C$ to a $C_l$ via a complete
intersection of two surfaces of degrees $f$ and $g$ satisfying
$\HH^1(\sI_C(v))=0$ for $v=f,g,f-4$ and $g-4$, and then apply $(a)$ to $C_l$.
To see that $C_l$ satisfies the assumption of $(a)$, first note that $M(C_l)$
admits a decomposition $M(C_l) \cong M'(C_l) \oplus M_{[f+g-4-t]}$ as
$R$-modules. Indeed $M=M(C)$ satisfies the duality 
\begin{equation} \label{prop41}
M(C_l) \cong {\Ext_R^4}(M,R)(-f-g) \cong \Hom_k(M,k)(-f-g+4) \ ,
\end{equation}
(cf. \cite{R} and \cite{MIG}, p.\! 133). If we let $M'(C_l) :=
{\Ext_R^4}(M',R)(-f-g)$, then the decomposition $M \cong M' \oplus M_{[t]}$
translates to
$$M(C_l) \cong {\Ext_R^4}(M',R)(-f-g) \oplus {\Ext_R^4}(M_{[t]},R)(-f-g) \cong
M'(C_l) \oplus M_{[t]}(2t+4-f-g) $$ since $ M_{[t]}(t) \cong
{\Ext_R^4}(M_{[t]}(t),R(-4))$ by the self-duality of the minimal resolution of
$M_{[t]}(t)$. Finally since $ M_{[t]}(2t+4-f-g)$ is supported in degree
$f+g-4-t$, we may write the module $ M_{[t]}(2t+4-f-g)$ as $
M(C_l)_{[f+g-4-t]}:= M_{[f+g-4-t]}$. Next to see that direct free part $F_1$
generated in degree $t$ in the resolution of $I(C)$, is equal (at least
dimensionally) to the corresponding part in degree $f+g-4-t$ of $F_2(C_l)(4)$
in the minimal resolution of $I(C_l)$ of the linked curve $C_l$, we remark
that since the isomorphism of \eqref{thm66} is given by the duality used in
\eqref{prop41}, it must commute with their decomposition as $R$-modules, i.e.
we have
\begin{equation} \label{prop42}
  {_0\!\Hom_R}(I(C), M(C)_{[t]}) \cong {_0\!\Hom_R}(M(C_l)_{[f+g-4-t]},E(C_l))
\end{equation} 
Then we conclude by \eqref{cor61} and \eqref{cor63} provided we can use
Remark~\ref{remmcor} for $v=0$. Indeed if $L^*:= \Hom_R(L,R)$, we have an
exact sequence $ \rightarrow (L_2')^* \rightarrow (L_3')^* \rightarrow
(L_4')^* \rightarrow {\Ext_R^4}(M',R) \cong M'(C_l)(f+g) \rightarrow 0$. Since
$L_2'$ has no direct free factor of degree $t$, it follows that
$(L_2')^*(-f-g)$ has no direct free factor of degree $f+g-t$, i.e. we have
${_{-4}\!\Hom_R}((L_2')^*(-f-g),M_{[f+g-4-t]})=0$ and Remark~\ref{remmcor}
applies. Now using $(a)$ to the linked curve $C_l$ with $m_2=m_1$, we get a
generization of $C_l'$ with constant postulation where $R(-f-g+t)^{m_1}$ is
"removed" in its minimal resolution. A further linkage, using a complete
intersection of the same type as in the linkage above (such a complete
intersection exists by \cite{K3}, Cor.\! 3.7) and the formula \eqref{prop42}
(replacing $C$ and $C_l$ by $C'$ and $C_l'$), gives the desired generization
$C'$, and we are done.
\end{proof}

\bigskip \bigskip

\end{document}